\documentclass[12pt,a4paper]{article}

\usepackage[english]{babel}
\usepackage{verbatim}
\usepackage{bm}

\usepackage{amsmath}
\usepackage{amssymb}
\usepackage{amsfonts}
\usepackage{amsbsy}

\usepackage{graphicx}
\usepackage{graphics}
\usepackage{epsfig}

\usepackage{array}
\usepackage{booktabs} 
\usepackage{caption}  
\usepackage{tabularx}
\usepackage{hyphenat} 
\usepackage{amsthm} 

\makeatletter
\def\@endtheorem{\endtrivlist}
\makeatother

\newtheorem{theorem}{Theorem}
\newtheorem{lemma}{Lemma}
\newtheorem{defin}{Definition}

\newtheorem*{proof*}{Proof}

\newtheorem*{remark*}{Remark}

\numberwithin{equation}{section}
\numberwithin{theorem}{section}
\numberwithin{lemma}{section}
\numberwithin{defin}{section}
\numberwithin{cor}{section}
\numberwithin{prop}{section}

\begin{document}

\title{Monotonicity, positivity and strong stability of the TR-BDF2  method and of its SSP extensions}
\author{Luca Bonaventura$^{(1)}$,\ \ Alessandro Della Rocca$^{(1),(2)}$}
\maketitle

\begin{center}
{\small
$^{(1)}$ MOX -- Modelling and Scientific Computing, \\
Dipartimento di Matematica, Politecnico di Milano \\
Via Bonardi 9, 20133 Milano, Italy\\
{\tt luca.bonaventura@polimi.it}\\
{\tt alessandro.dellarocca@polimi.it}
}
\end{center}

\begin{center}
{\small
$^{(2)}$ Tenova S.p.A., \\
Global R\&D, \\
Via Albareto 31, 16153 Genova, Italy\\
{\tt alessandro.dellarocca@tenova.com}
}
\end{center}

\date{}


\vspace*{1.0cm}


\vspace*{1.0cm}

\pagebreak

\abstract{We analyze the one-step method TR-BDF2 from the point of view of monotonicity, strong stability and
  positivity. All these properties are strongly related and reviewed in the common framework of absolute
  monotonicity. The radius of absolute monotonicity is computed and it is shown that the parameter value which
  makes the method L-stable is also the value which maximizes the radius of monotonicity. Two hybrid variants
  of TR-BDF2 are proposed, that reduce the formal order of accuracy and maximize the absolute monotonicity
  radius, while keeping the native L-stability useful in stiff problems. Numerical experiments compare these
  different hybridization strategies with other methods commonly used in the presence of stiff and mildly
  stiff source terms. The results show that both strategies provide a good compromise between accuracy and robustness
  at high CFL numbers, without suffering from the limitations of alternative approaches already available in
  literature.}

\pagebreak

\section{Introduction}
 \label{intro} \indent

The classical error analysis of numerical methods for ordinary differential equations (ODE) does not yield
sufficient criteria for the conservation of special properties of the continuous solutions during time
integrations. In many applications, for instance, the solution is required to remain non negative or to take
values in a certain range or again to preserve the total variation as a function of the space variable, in
case a time and space dependent partial differential equation is solved numerically.  A number of
different strategies have been proposed to address each of these issues separately, see
e.g. \cite{burchard:2003}, \cite{hundsdorfer:2003a}, \cite{sandu:2001}, \cite{shu:1988a}.  Especially in the
ODE literature, many of these problems are empirically resolved by step size adaptation strategies that complement
traditional ODE methods. However, in many applications, such as numerical weather prediction, environmental
fluid dynamics or turbulent reactive flow simulations, the step size is usually kept fixed and/or relatively
large, in order to minimize the number of expensive function evaluations required from the complex source
terms involved.  While dynamic time step adaptation \cite{hairer:1996} and multirate approaches
\cite{constantinescu:2007} are able to overcome these problems, in this work we will study a robust but
accurate fixed time step approach that can guarantee a good compromise between preservation of some relevant
monotonicity properties, accuracy and efficiency.  Splitting approaches \cite{sportisse:2000} are commonly used to
couple complex source terms to the discretized fluid dynamics in a relatively simple way, but we will restrict
our attention to methods that do not resort to operator splitting, which may entail a loss of accuracy for
advection-diffusion-reaction problems with space dependent source terms.
 
We will focus on the analysis of the monotonicity properties of the TR-BDF2 method, that was introduced in
\cite{bank:1985} and successively reformulated and analyzed in \cite{hosea:1996}. This second order accurate,
A-stable and L-stable method has a number of interesting properties and it has been recently used with success
in \cite{tumolo:2015} as the key ingredient of an efficient semi-implicit, semi-Lagrangian discretization of
fluid dynamics equations representative of many environmental models. It is therefore of interest to
understand to which extent this method can also guarantee positivity and monotonicity for the equations of
advected species. We will show that TR-BDF2 is conditionally monotone under a time step
restriction that allows for time steps more than double with respect to those of explicit schemes. Our analysis
relies on the results in \cite{ferracina:2004}, \cite{ferracina:2005}, \cite{gottlieb:2001},
\cite{higueras:2004}, \cite{kraaijevanger:1991}, \cite{spijker:2007}, that allow to derive sufficient
conditions for properties like positivity, monotonicity and total variation preservation in a unified
framework of an extended concept of monotonicity. We then propose two different modifications of the TR-BDF2 method,
both based on a hybridization with the unconditionally monotone implicit Euler method. In this way, accuracy
is sacrificed locally in space or time in order to preserve monotonicity, independently of the time step and
stiffness of the problem. Other approaches, focusing specifically on
the equations of chemical kinetics, have been proposed under more
restrictive assumptions in \cite{bertolazzi:1996}, \cite{formaggia:2011}, \cite{vanveldhuizen:2006}. 
The present approach represents an improvement over these results,  since it does not require  source term splitting,
it is not limited to non stiff problems as \cite{bertolazzi:1996}, nor it requires a special form of the source
term as \cite{formaggia:2011} and differently than \cite{vanveldhuizen:2006} it is second order when critical
solution properties are not violated under the selected time step size.
   
In Section~\ref{review}, the theory of monotone and SSP methods is reviewed. The TR-BDF2 method is presented
in detail in Section~\ref{trbdf2} where its absolute monotonicity property is analyzed. Two strategies to
improve its monotonicity properties irrespective of time step size are described in
Section~\ref{trbdf2_variants} stemming from the classical results of absolute monotonicity. Other competitive
time discretization approaches are introduced in Section~\ref{competitors} and interpreted under the SSP
theory whenever possible.  In Section~\ref{experiments} we present an empirical assessment of the properties
of all these time integration methods in a variety of relevant benchmark problems. Conclusions, results and
directions of further investigations are summarized in Section~\ref{conclusions}.

\section{Review of monotonicity and strong stability results}
 \label{review}

We review the recent progresses in the field of strong stability preserving (SSP) methods introduced, among others, in
\cite{ferracina:2004}, \cite{gottlieb:2001}, \cite{gottlieb:2011}, \cite{higueras:2004} and \cite{shu:1988b}.
  In this work we consider an initial value problem for a system of ordinary differential
equations (ODEs) of type
\begin{subequations}
 \label{eq:ode}
  \begin{align}
   u^{\prime}(t)&=f(t,u(t)) \quad \text{and} \quad t\in[0,T], \\
   u(0) &= u^0.
  \end{align}
\end{subequations}
We assume that   $u^0 \in \mathbb{R}^m$, $f : \mathbb{R}\times\mathbb{R}^m \to \mathbb{R}^m$
such that the problem \eqref{eq:ode} has a unique solution. Moreover we assume also that
$\|\!\cdot\! \| : \mathbb{R}^m \to \mathbb{R}$ is a convex functional
 $$\| \lambda v + \left( 1 - \lambda \right) w \| \leq \lambda \| v \| + \left( 1 - \lambda \right) \| w \| $$
for $0 \leq \lambda \leq 1 $ and $v$,$w \in \mathbb{R}^m$.
We shall deal with numerical methods for finding a numerical approximation $u^n$ to the exact solution values
$u(n \, h)$, where $h$ is a positive step size, assumed for simplicity to be constant and equal to
$h = T/N_t,$ with $T$ the final time of integration and $n=1,\dots,N_t$.
 
Monotonicity of the total variation seminorm, discrete maximum principle, positivity or range boundedness are
all nonlinear properties that can be seen as specific realizations of some form of monotonicity. Even
contractivity, a property relevant for the numerical stability of time integration methods, can be
reintepreted as a form of monotonicity. The recent theory of SSP integration methods \cite{gottlieb:2011}
shows that all these concepts are strongly related. As a consequence, they can be guaranteed during the
numerical integration using the same fundamental approach. In this section we provide a synthetic presentation
of this framework, that is based on the classical results in \cite{kraaijevanger:1991}. We will focus in
particular on the properties of $s$-stages Runge-Kutta methods (RK)
\begin{subequations}
 \label{eq:RK}
\begin{align}
   g^i &= u^{n} + h \sum\limits_{j=1}^s a_{ij}\,f(t^n + c_j h,g^j)\quad(i=1,\ldots,s) \\
   u^{n+1} &= u^{n} + h \sum\limits_{i=1}^s b_i\,f(t^n + c_i h,g^i) 
\end{align}
\end{subequations}
where $a_{ij}$, $b_i$ and $c_i$ are real parameters which characterize the method and $g^i$ are the intermediate
stages. It is usually assumed that $\sum\limits_{j=1}^s a_{ij} = c_i$.  The method is explicit if
$a_{ij}=0$ for $j \geq i$ and implicit otherwise. The parameters of the method are traditionally collected in
compact form in the Butcher tableau as an $s \times s$ matrix $A = (a_{ij})$, a row vector $b =
(b_1,\ldots,b_s)^\intercal$ and a column vector $c = (c_1, \ldots,c_s)$. For homogeneous initial value
problems, the coefficients $c_i$ are not relevant, so that each RK method is completely identified by its
coefficients $(A,b)$.

\begin{defin}
\label{def:discrete_monotonicity}
  \textbf{Monotonicity}. The RK method \eqref{eq:RK} is monotone with respect to the functional
   $\|\!\cdot\! \|$ if $ \|u^{n}\| \leq \|u^{0}\|$ 
   under the assumption that
   \begin{equation}
    \label{eq:hp_discrete_monotonicity}
     \|u + h f(t,u)\| \leq \|u\| \quad\text{for}\quad 0 < h \leq \tau_{0}.
   \end{equation}
\end{defin}
Related definitions also
commonly found in the literature   are internal monotonicity
\begin{equation}
 \label{eq:internal_monotonicity}
   \|g^{i}\| \leq \|u^{n}\| \quad\text{for}\quad 1 \leq i \leq s
\end{equation}
and external monotonicity 
\begin{equation}
 \label{eq:external_monotonicity}
  \|u^{n+1}\| \leq \|u^{n}\|.
\end{equation}
Assumption \eqref{eq:hp_discrete_monotonicity}, commonly found in many references, see e.g.
\cite{ferracina:2004}, \cite{higueras:2004}, \cite{higueras:2005a}, \cite{higueras:2013},
\cite{hundsdorfer:2003b}, \cite{hundsdorfer:2011a}, \cite{ketcheson:2011}, \cite{shu:1988a} and
\cite{shu:1988b}, essentially amounts to define $\tau_0$ as the maximum time step under which the explicit
Euler method is monotone. In these studies the critical step size for monotonicity is determined such that
property \ref{def:discrete_monotonicity} is verified for all
\begin{equation}
 \label{eq:monotone_stepsize}
  0 < h \leq c \, \tau_0
\end{equation}
thus making the RK method \emph{conditionally monotone}. If property \ref{def:discrete_monotonicity} is
verified for any step size $h,$ then the method is called \emph{unconditionally monotone}. In assessing
monotonicity of different RK methods, the interest lies usually in determining the
\emph{maximal step size coefficient} $c$ such that a time integration method is conditionally monotone.
 
Frequently, the convex functional $\| \!\cdot\! \|$ is intended either as the  supremum norm $\| x \| =
\|x\|_{\infty} = \sup_i | x_i |$ or as the total variation seminorm $\| x \| = \|x\|_{TV} = \sum_i |
x_{i+1} - x_{i} |, $ where $x_i$ are the components of the vector $x$. We remind that $\|x\|_{TV}$ is a
seminorm since it may vanish for $x \neq 0$ when $x_{i} = C$.
Numerical methods statisfying \ref{def:discrete_monotonicity} under the total variation seminorm are
called \emph{total variation diminishing (TVD)}. They are especially important in the numerical solution of
hyperbolic conservation laws, see e.g.  \cite{hundsdorfer:2003a}, \cite{leveque:2002}, \cite{shu:1988b}.

For some classes of initial value problems, the properties of positivity and range boundedness play important
roles in obtaining physically meaningful numerical solutions. Furthermore, due to the strongly nonlinear form
of such problems, the ability in maintaining such native properties also in the numerical solutions is
important in order to guarantee numerical stability of time integrations.

\begin{defin}
 \label{def:positivity}
  \textbf{Positivity}. The RK method \eqref{eq:RK} is positive if whenever $u^{0} \geq 0$
  it guarantees that $u^{n+1} \geq 0$
  under the assumption that
   \begin{equation}
    \label{eq:hp_positivity}
     u + h f(t,u) \geq 0 \quad\text{for}\quad 0 < h \leq \tau_{0}.
   \end{equation}
\end{defin}

\begin{defin}
 \label{def:range_boundedness}
  \textbf{Range boundedness}. The RK method \eqref{eq:RK} is range bounded in $[\chi,\psi]$ if
  whenever $\chi \leq u^{0} \leq \psi$ it guarantees that $\chi \leq u^{n+1} \leq \psi$ 
  under the assumption that
   \begin{equation}
    \label{eq:hp_range_boundedness}
     \chi \leq u + h f(t,u) \leq \psi \quad\text{for}\quad 0 < h \leq \tau_{0}.
   \end{equation}
\end{defin}
Both these properties are usually guaranteed if a  time step limitation of the form
\eqref{eq:monotone_stepsize} is respected.  Even though these properties are formally different from the monotonicity
property \ref{def:discrete_monotonicity}, they can be equally derived from monotonicity after proper
assumptions on the function $f$. In this respect, the generalization of monotonicity to arbitrary sublinear
functionals $\|\!\cdot\! \|$ becomes relevant.  Following the presentation in \cite{hundsdorfer:2011b} and
\cite{spijker:2007}, it is useful to introduce two sublinear functional, denoted  respectively as
\emph{floor} and \emph{ceil} functional
\begin{subequations}
 \label{eq:functionals_floor_ceil}
  \begin{align}
   \| u \|_{\lfloor \chi \rfloor} &= -\min_j (\chi, u_j) \\
   \| u \|_{\lceil \psi \rceil} &= \max_j (\psi, u_j).
  \end{align}
\end{subequations}
These functionals are not seminorms, since they both violate property $\|\lambda v\| = |\lambda| \|v\|$ for
$\lambda\!=\!-1$. Using both functionals the range boundedness property \ref{def:range_boundedness} naturally
follows, while setting $\chi \!=\! 0$ in the floor functional the positivity property \ref{def:positivity} is
recovered.  As a consequence, by introducing the floor and ceil functionals it possible to recast Definitions
\ref{def:positivity} and \ref{def:range_boundedness} in a form similar to Definition
\ref{def:discrete_monotonicity}. Thus positivity and range boundedness can just be interpreted as different
forms of monotonicity.

Alternatively, positivity and range boundedness can also be considered as two alternative realizations of the
discrete maximum principle.
\begin{defin}
 \label{def:discrete_maximum_principle}
  \textbf{Discrete Maximum Principle}. The RK method \eqref{eq:RK} follows the discrete maximum principle if
  it guarantees that 
  \begin{equation*}
   \min_j u_{j}^{0} \leq u_{j}^{n+1} \leq \max_j u_{j}^{0}
  \end{equation*}
  under the assumption that for $0 < h \leq \tau_{0}$ and
  $\forall u \in \mathbb{R}^m$ with components $u_p$ 
  \begin{equation}
   \label{eq:hp_discrete_maximum_principle}
    \min_{1 \leq q \leq m} u_q \leq u_p + \tau_0 f_p(t,u(t)) \leq \max_{1 \leq q \leq m} u_q, \quad (1 \leq p \leq m).
  \end{equation}
\end{defin}
Similarly to the other properties, range boundedness may be verified under a step size restriction analogous to
 \eqref{eq:monotone_stepsize}.  Again following \cite{hundsdorfer:2011b} and \cite{spijker:2007}, we introduce
 two relevant sublinear functionals, denoted as \emph{max} and \emph{min} functional
\begin{subequations}
 \label{eq:functionals_+-}
  \begin{align}
   \| u \|_> &= \max_j u_j \\
   \| u \|_< &= -\min_j u_j 
  \end{align}
\end{subequations}
which allow us to write the assumption \eqref{eq:hp_discrete_maximum_principle} in the form
\eqref{eq:hp_discrete_monotonicity}. By assuming the monotonicity property \ref{def:discrete_monotonicity}
under the functionals \eqref{eq:functionals_+-} the discrete maximum principle
\ref{def:discrete_maximum_principle} directly follows.

Contractivity of numerical approximation methods has been extensively studied. Relevant conclusions on
step size conditions for contractivity have been given in \cite{spijker:1983}, while contractivity of RK for
nonlinear problems was thoroughly examined in \cite{kraaijevanger:1991}.

\begin{defin}
 \label{def:contractivity}
  \textbf{Contractivity}. The RK method \eqref{eq:RK} is contractive if
   $\|\tilde u^{n+1} - u^{n+1}\| \leq \|\tilde u^{0} - u^{0}\|$
  under the assumption that
   \begin{equation}
    \label{eq:hp_contractivity}
     \|\tilde u - u + h (f(t,\tilde u) - f(t,u))\| \leq \|\tilde u - u\| \quad\text{for}\quad 0 < h \leq \tau_{0}.
   \end{equation}
\end{defin}
Usually, Definition \ref{def:contractivity} is verified under a step size restriction in the form of
\eqref{eq:monotone_stepsize}. For conditional contractivity, the \emph{circle condition} that was originally
assumed in \cite{kraaijevanger:1991} is
\begin{equation}
 \label{eq:circle_condition}
  \|f(t,\tilde u) - f(t,u) + \rho(\tilde u - u) \| \leq \rho \|\tilde u - u\|
\end{equation}
It was shown later in \cite{higueras:2005a} that this condition can be considered as a special form of
\eqref{eq:hp_discrete_monotonicity}.
%
See also \cite{spijker:2007} for additional considerations on these issues.

In the framework introduced above, all these properties appear as different forms of monotonicity under a
proper choice of the convex functional $\|.\|$. It is thus possible to extend the analytical results from the
preservation of a specific property to that of any other related property. Every convex functional will be retained under convex
combinations of single RK methods preserving it.  Even though in principle this approach may not lead to sharp
bounds, experience shows that the time step limits obtained from monotonicity are representative for the
preservation of the other properties as well \cite{hundsdorfer:2011a}. The sharper results obtained in the
context of inner product norms \cite{higueras:2005b}, for example, did not lead to time step limits
significantly different for the practical use of conditionally monotone RK methods. It is to be noticed that
  RK methods that are non-monotone under arbitrary functionals may indeed be conditionally monotone under inner product norms, as in
\cite{higueras:2005b}.  For linear problems these conditions relax even further, see e.g., \cite{spijker:1983},
\cite{vandegriend:1986}.

As a consequence, it appears useful to extend the contractivity results to the other properties. To this
purpose, the absolute monotonicity introduced in \cite{kraaijevanger:1991} allows to investigate necessary and
sufficient conditions for contractivity of RKs in dissipative problems under sublinear functionals, including
the maximum norm. In \cite{ferracina:2004} this approach was extended to study monotonicity of general RK
under arbitrary seminorms. One of the relevant results is that the maximal step size coefficient for
monotonicity is equal to the maximal step size for contractivity.  Recently these conclusions were extended to
the analysis of general linear methods under arbitrary convex functionals in \cite{spijker:2007}.

%

For our next discussion we will focus on \emph{irreducible} RK methods, which are the only ones practically
relevant. For a definition of irreducibility see \cite{hairer:1993}.
Following  \cite{kraaijevanger:1991}, we introduce for real $\xi$ the quantities
\begin{equation}
 \label{eq:a_b_e_phi}
  \begin{aligned}
   A(\xi) &= A(I - \xi A)^{-1}, & \quad b^\intercal(\xi) &= b^\intercal(I - \xi A)^{-1}, \\
   e(\xi) &= (I - \xi A)^{-1}e, & \quad \varphi(\xi) &= 1 + \xi b^\intercal(I - \xi)A^{-1}e.
  \end{aligned}
\end{equation}

\begin{defin}
 \label{def:absolute_monotonicity}
  \textbf{Absolute monotonicity of RK}. An irreducible $s$-stage RK   $(A,b)$ is absolutely
  monotone at $\xi \in \mathbb{R}$ if $A \geq 0$, $b \geq 0$, $e \geq 0$ and $ \varphi \geq 0$ elementwise.
\end{defin}
For the stability function $\varphi,$ this entails that $\frac{\mathrm{d}^k \varphi}{\mathrm{d}z^k}(\xi) \geq 0$
for any $k \geq 0$, since the rationale behind Definition \ref{def:absolute_monotonicity} lies in the Taylor
expansion of some characteristic functions of a RK method, including the stability function. The quantities
\eqref{eq:a_b_e_phi} form the coefficients of such expansions, see \cite{higueras:2006a} and
\cite{kraaijevanger:1991} for further details. These notations are useful to introduce the radius of absolute monotonicity for a RK method.
\begin{defin}
 \label{def:radius_absolute_monotonicity}
  \textbf{Radius of absolute monotonicity}. An $s$-stage RK with scheme $(A,b)$ and $A \geq 0$ and $b \geq 0$
  is characterized by its radius of absolute monotonicity defined for all $\xi$ in
  $-r \leq \xi \leq 0$ as
  \begin{equation}
   \label{eq:radius_abs_monotonicity}
    \begin{split}
     R(A,b) = \sup\{r : &r \geq 0, \\
                      &A(\xi) \geq 0, \;  b^\intercal(\xi) \geq 0, \; e(\xi) \geq 0, \; \varphi(\xi) \geq 0\}.
          \end{split}
  \end{equation}
\end{defin}
In \cite{kraaijevanger:1991}, two useful results are derived, that simplify the practical estimation of the
absolute monotonicity radius. We introduce the incidence matrix $\text{Inc(A)}$ as the matrix containing $0$
or $1$ if the corresponding element in the matrix $A$ is $a_{i,j}=0$ or $a_{i,j} \neq 0$ respectively.

\begin{theorem}[\cite{kraaijevanger:1991}, Theorem 4.2]
 \label{monotonicity_conditions}
 For an irreducible RK $R(A,b)>0$ iff $A\geq0$, $b>0$ and $\text{Inc}(A^2) \leq \text{Inc}(A)$.
\end{theorem}

\begin{lemma}[\cite{kraaijevanger:1991}, Lemma 4.4]
  For an irreducible RK $R(A,b) \geq r$ iff $A \geq 0$ and $(A,b)$ is absolutely monotone at $\xi = -r$.
\end{lemma}

More compact definitions of the radius of absolute monotocity involve the use of matrices derived from the
Butcher tableau, see e.g., \cite{ketcheson:2009}, but for the purpose of this work we found the use of the
original definition more convenient.  Other classical results for the practical identification of RK
properties from the corresponding Butcher tableau stem from the \emph{stage order}, that has practical
relevance to avoid the order reduction phenomenon during stiff transients.

\begin{theorem}[\cite{kraaijevanger:1991}, Theorem 8.5]
  Any RK with $A \geq 0$ has stage order $\tilde p \leq 2$. If $\tilde p = 2$ then $A$ must have a zero row.
\end{theorem}

\begin{lemma}[\cite{kraaijevanger:1991}, Lemma 8.6]
  Any RK with $b > 0$ has stage order $\tilde p \geq \lfloor \frac{p-1}{2}\rfloor\ $.
\end{lemma}
 
\begin{theorem}[\cite{kraaijevanger:1991}, Theorem 8.3]
  \textbf{Unconditional contractivity}.\\ The order barrier for $R(A,b) = \infty$ in a RK is $p \leq 1$. Some
  first order unconditionally contractive RKs are the implicit Euler, the RADAU IA and the RADAU IIA methods.
\end{theorem}

\begin{theorem}[\cite{kraaijevanger:1991}, Corollary 8.7]
  \textbf{Conditional contractivity}.\\ The order barriers for $R(A,b) \geq 0$ under the circle condition
  \eqref{eq:circle_condition} are $p \leq 4 $  for explicit RK  and  $p \leq 6$
  for implicit RK.
\end{theorem}
As a consequence, step size restrictions on RK of formal order $p > 1$ are inevitable to preserve contractivity
of the numerical solution.  The analysis was recently extended in \cite{spijker:2007} to a much larger class of
methods, including IMEX methods (see also \cite{higueras:2006a}). This implies that the order barriers above
are inevitable for a very large class of time discretization methods.

In recent literature \cite{gottlieb:2011}, \cite{gottlieb:2001}, methods satisfying condition
\eqref{eq:external_monotonicity} for a general convex functional are called \emph{strong\hyp{}stability
  preserving (SSP)}, in order to specify their ability to preserve any convex functional bound. Thus they
generalize classical TVD methods specifically developed for hyperbolic conservation laws.

\begin{defin}
\label{def:SSP}
  \textbf{Strong stability preserving (SSP)}. The RK method \eqref{eq:RK} is SSP with respect to the
  functional $\|\!\cdot\! \|$ if $ \|u^{n}\| \leq \|u^{0}\|$ under the assumption that
  \begin{equation}
   \label{eq:hp_SSP}
    \|u + h f(t,u)\| \leq \|u\| \quad\text{for}\quad 0 < h \leq \tau_{0}.
  \end{equation}
  The \emph{SSP coefficient} is the largest constant $c \geq 0$ such that this definition is verified
  for all $0 < h \leq c \tau_{0}$.
\end{defin}
The definition above closely follows Definition \ref{def:discrete_monotonicity} for monotonicity; indeed,
these definitions  are  equivalent.
%
%
%
The SSP coefficient $c$ turns out to be strongly related to the radius of absolute monotonicity introduced in
Definition \ref{def:radius_absolute_monotonicity}.
\begin{theorem}[\cite{ferracina:2005}, Theorem 3.4]
  For an irreducible RK, $c = R(A,b)$.
\end{theorem}
Recent studies focused on the search of SSP-optimal RK methods having large SSP coefficients. While explicit
SSP RKs are known since the seminal work in \cite{shu:1988b}, a search for implicit SSP RKs started only
recently in \cite{ferracina:2008} and \cite{ketcheson:2009}, where it was found that the SSP-optimal implicit
RKs of order $p=2$ and $p=3$ are indeed SDIRK, while the optimal methods for $p=4$ are DIRK. In the quest for
improved monotonicity, the classes of \emph{two-step RKs} \cite{ketcheson:2011} and \emph{diagonally split RK}
\cite{bellen:1994}, \cite{bellen:1997}, \cite{macdonald:2008} have been investigated in the literature, with mixed success.
Starting from the framework introduced above in our following analysis we will consider all these properties
just as specific realizations of absolute monotonicity.  As a consequence, in the next sections they will
be briefly referred to as \emph{monotonicity} of the time integration methods.

\section{Strong stability preservation for the TR-BDF2  method}
 \label{trbdf2}

We analyze here the monotonicity properties of the TR-BDF2 method, originally introduced in \cite{bank:1985}
and successively studied in \cite{hosea:1996} as a DIRK. The same method was rediscovered in
\cite{butcher:2000} and has been applied also in \cite{giraldo:2013} to treat the implicit terms in an
additive Runge-Kutta (ARK).  A semi-implicit, semi-Lagrangian reinterpretation of this method has recently
been proposed  in \cite{tumolo:2015}.
 
We rely on the monotonicity and contractivity results reviewed in Section~\ref{review}.  In its original
formulation, TR-BDF2 is defined as a one-step method resulting from the composition of the trapezoidal rule in
the first substep, followed by BDF2 in the second substep. This combination is empirically justified under the
rationale of combining the good accuracy of the trapezoidal rule with the stability and damping of fast modes
guaranteed by BDF2. The TR-BDF2 method is
\begin{subequations}
 \label{eq:trbdf2_original}
  \begin{align}
   u^{n+\gamma} - \frac{\gamma}2 h f^{n+\gamma} &= u^n + \frac{\gamma}2 h f^n \\
   u^{n+1} - \frac{(1-\gamma)}{(2-\gamma)} h f^{n+1} &= \frac{1}{\gamma(2-\gamma)}u^{n+\gamma} -
      \frac{(1-\gamma)^2}{\gamma(2-\gamma)} u^n
  \end{align}
\end{subequations}
where $\gamma \in (0,1) $ is a parameter whose value determines the stability and monotonicity properties of the
method. By requiring that both stages have the same Jacobian, the value $\gamma = 2 - \sqrt 2$
was derived in \cite{bank:1985}. This value is also the only one for which the method \eqref{eq:trbdf2_original}
is L-stable. As outlined in \cite{hosea:1996}, the TR-BDF2 method can be rewritten as a DIRK method. However,
instead of closely following this reference, we will first derive the DIRK family associated to
\eqref{eq:trbdf2_original} before imposing the condition on the Jacobian.
The Butcher tableau for this DIRK family is then
\begin{center}
\begin{tabular}{c|ccc}
  $0$ & $0$ & & \\ [0.5ex]
  $\gamma$ & $\frac{\gamma}{2}$ & $\frac{\gamma}{2}$ & \\ [1.0ex]
  $1$ & $\frac{1}{2(2-\gamma)}$ & $\frac{1}{2(2-\gamma)}$ & $\frac{1-\gamma}{2-\gamma}$ \\ [1ex]
  \hline \\ [-1.5ex]
  & $\frac{1}{2(2-\gamma)}$ & $\frac{1}{2(2-\gamma)}$ & $\frac{1-\gamma}{2-\gamma}$
\end{tabular}
\end{center}
We remark that, by imposing the condition $a_{22} = a_{33},$ the optimal value for $\gamma$ able to guarantee
L-stability is readily obtained. The stability function associated to the above DIRK family is
\begin{equation}
 \label{eq:trbdf2_stability_function}
 \varphi(\xi) =\frac{\det (I - \xi A + \xi e b^\intercal)}{\det (I - \xi A)}  
 =\frac{[1 + (1 - \gamma)^2] \xi + 2(2-\gamma)}
      {2(2 - \gamma)(1 - \xi \frac{\gamma}{2})(1 - \xi \frac{1-\gamma}{2-\gamma})}.
\end{equation}
The rational polynomial defining the stability function of the DIRK family has a single zero at $\xi \!=\!
\frac{2(\gamma-2)}{1+(1-\gamma)^2}$, which is always negative for all possible values of parameter $\gamma$,
and two poles at $\xi \!=\! \frac{2}{\gamma}$ and $\xi \!=\! \frac{2-\gamma}{1-\gamma}$ respectively, which
are always positive. While in \cite{bank:1985} it is argued that there is a double pole at $\xi \!=\!
\frac{2}{\gamma}$, this is true only for the numerical value of $\gamma$ satisfying $a_{22}\!=\!a_{33}$.  As
shown in \cite{hosea:1996}, the TR-BDF2 method is embedded in a (2,3) Runge-Kutta pair, thus allowing an
efficient estimation of the time discretization error, in case  step size adaptation is required.  Additionally, it is
immediate to find out that the stage order of TR-BDF2 is indeed $\tilde p \!=\! 2$, thus making the method
resilient to order reduction in stiff problems (see e.g., \cite{higueras:2004}, \cite{macdonald:2008}).

We analyze the monotonicity properties of the DIRK family generalizing the TR-BDF2 by following
\cite{ferracina:2004}, \cite{ferracina:2005}, \cite{kraaijevanger:1991}. Following the definition
\eqref{eq:radius_abs_monotonicity}, we find

\begin{equation*}
 \label{eq:a_trbdf2}
A(\xi)=
\begin{bmatrix}
 0 & 0 & 0
 \\[0.5em]
 \frac{\frac{\gamma}{2}}{1-\frac{\gamma}{2}\xi} & \frac{\frac{\gamma}{2}}{1-\frac{\gamma}{2}\xi} & 0
 \\[0.5em]
 \frac{1}{2(2-\gamma)(1-\frac{\gamma}{2}\xi)\beta} 
 & \frac{1}{2(2-\gamma)(1-\frac{\gamma}{2}\xi)\beta} 
 & \frac{1-\gamma}{(2-\gamma)\beta}
\end{bmatrix}
\end{equation*}

\begin{equation*}
 \label{eq:b_trbdf2}
b^\intercal(\xi)=
\begin{bmatrix}
 \frac{\beta -\xi\gamma/2}{2(2-\gamma)}  
 \\[0.5em]
 \frac{\beta -\xi\gamma/2}{2(2-\gamma)}  
 \\[0.5em]
 \frac{(1-\gamma)\beta}{2-\gamma}  
\end{bmatrix}
%
\ \ \ \ \ 
e(\xi)=
\begin{bmatrix}
 1
 \\[0.5em]
 \frac{1+\frac{\gamma}{2}\xi}{1-\frac{\gamma}{2}\xi}
 \\[0.5em]
 \frac{[1 + (1 - \gamma)^2] \xi + 2(2-\gamma)}
      {2(2 - \gamma)(1 - \xi \frac{\gamma}{2})\beta}
\end{bmatrix}
\end{equation*}

\begin{equation*}
 \label{eq:phi_trbdf2}
\varphi(\xi)=
 \frac{[1 + (1 - \gamma)^2] \xi + 2(2-\gamma)}
      {2(2 - \gamma)(1 - \xi \frac{\gamma}{2})\beta}
\end{equation*}
where we have set $ \beta \!=\! 1-\xi(1-\gamma)/(2-\gamma).$ The main conclusion is that the radius of
absolute monotonicity of TR-BDF2 is $R(A,b) \!=\! \frac{2(2-\gamma)}{1+(1-\gamma)^2}$ and that it is maximized
for $\gamma \!=\! 2-\sqrt{2}$, which is exactly the value that makes the DIRK family L-stable. With this value
of the parameter $\gamma,$ the absolute monotonicity radius is $R(A,b) \!\approx\! 2.414$.  From
\cite{ferracina:2004}, we conclude that this is also the step size coefficient $c$ for conditional monotonicity
of the method under arbitrary seminorms and sublinear functionals for any nonlinear problem.

\section{Two unconditionally monotone variants of TR-BDF2}
 \label{trbdf2_variants}

The results reviewed in Section~\ref{review} lead to the conclusion that there are no unconditionally monotone
RK methods of order higher than one.  For stiff initial value problems, this implies that even implicit higher order RK 
methods will
be always subject to a CFL-like condition for monotonicity.  Thus, the only RK method  that does not need to
comply with a time step restriction and that can be safely used without time step adaption is in practice  the implicit
Euler method. Due to such limitation, which is particularly relevant for problems with chemical kinetics, some
solution methods found in literature rely exclusively on it, sacrificing accuracy for improving stability and
consistency \cite{vanveldhuizen:2008}.

We propose here two hybridization strategies of TR-BDF2 with the implicit Euler method, that can be activated
using a sensor detecting violations of relevant functional bounds.
These hybrid schemes bring time integration back to a first order unconditionally monotone method whenever
the sensor detects a violation of a selected functional bound during the current integration. Empirically, this
local loss in accuracy should not be too detrimental on solution accuracy if time step is not too large.  Our
proposed methods apply this safe mode whenever it is expected   from SSP theory that TR-BDF2 may produce
non-monotone solutions.

The \textbf{hybrid TR-BDF2} method is thus obtained by introducing a weighting parameter $\alpha \in \, [0,1]$
in both   stages of TR-BDF2
\begin{subequations}
 \label{eq:trbdf2_hybrid}
  \begin{align}
   u^{n+\gamma} -\gamma h (1-\frac{\alpha}{2}) f^{n+\gamma}  &=  u^n+ \gamma h \frac{\alpha}{2} f^n \\
   u^{n+1} -\frac{(1-\gamma) h}{\alpha(1-\gamma)+1} f^{n+1} &=  \nonumber \\
    &\kern -1.2in =\frac{\alpha(\frac{1}{\gamma}-1)+1}{\alpha(1-\gamma)+1}  u^{n+\gamma} 
        -\frac{\alpha}{\alpha(1-\gamma)+1}  \frac{(1-\gamma)^2}{\gamma} u^{n}.
  \end{align}
\end{subequations}
%
%
%
%
%
For $\alpha=1$ the first step is the trapezoidal rule and the second step is the BDF2 formula, thus
reconstructing the original TR-BDF2 \eqref{eq:trbdf2_original}. For $\alpha=0$ each of the two steps above is
equivalent to an implicit Euler step, thus transforming the hybrid TR-BDF2 in succession of two substeps of
the implicit Euler method (IE-IE)  of length $\gamma h$ and $(1- \gamma)h,$ respectively, and making the method
unconditionally monotone.

The hybrid TR-BDF2 method can be rewritten as a DIRK scheme as done for \eqref{eq:trbdf2_original}. By
injecting the first step in the second one, the Butcher tableau of the hybrid TR-BDF2 method is found
\begin{center}
\begin{tabular}{c|ccc}
  $0$ & $0$ & & \\ [0.5ex]
  $\gamma$ & $\gamma \frac{\alpha}{2}$ & $\gamma\left(1-\frac{\alpha}{2}\right)$ & \\ [1ex]
  $1$ & $\frac{\alpha}{2} \frac{\alpha(1-\gamma)+\gamma}{\alpha(1-\gamma)+1}$ & $\left(1 - \frac{\alpha}{2}\right) \frac{\alpha(1-\gamma)+\gamma}{\alpha(1-\gamma)+1}$ & $\frac{1-\gamma}{\alpha(1-\gamma)+1}$ \\ [1ex]
  \hline \\ [-1.5ex]
  & $\frac{\alpha}{2} \frac{\alpha(1-\gamma)+\gamma}{\alpha(1-\gamma)+1}$ & $\left(1 - \frac{\alpha}{2}\right) \frac{\alpha(1-\gamma)+\gamma}{\alpha(1-\gamma)+1}$ & $\frac{1-\gamma}{\alpha(1-\gamma)+1}$
\end{tabular}
\end{center}
For the unconditionally monotone method ($\alpha=0$), equivalent to a double step of implicit Euler, the
above above reduces to
\begin{center}
\begin{tabular}{c|ccc}
  $0$ & $0$ & & \\ [0.5ex]
  $\gamma$ & $0$ & $\gamma$ & \\ [0.5ex]
  $1$ & $0$ & $\gamma$ & $1-\gamma$ \\ [1ex]
  \hline \\ [-1.5ex]
  & $0$ & $\gamma$ & $1-\gamma$
\end{tabular}
\end{center}
The stability function of the DIRK family associated to the hybrid TR-BDF2 method is thus
\begin{equation}
\label{eq:trbdf2-blended_stability_function}
 \varphi(\xi) = \frac{ 1 + \left[ \frac{\alpha(1-\gamma)+\gamma}{\alpha(1-\gamma)+1}
                                 - \gamma\left(1-\frac{\alpha}{2} \right) \right] \xi }
    { 1 - \left[ \frac{1-\gamma}{\alpha(1-\gamma)+1}
                                 + \gamma\left(1-\frac{\alpha}{2}\right) \right] \xi
       + \gamma \left(1-\frac{\alpha}{2}\right) \frac{1-\gamma}{\alpha(1-\gamma)+1} \xi^2 }
\end{equation}
From this expression for the case $\alpha=1$ we recover the stability function
\eqref{eq:trbdf2_stability_function} and for the unconditionally monotone case $\alpha=0$ the stability
function $ \varphi(\xi) = 1/[1 - \xi + \gamma(1-\gamma) \xi^2] $.

\begin{figure}[htbp]
    \includegraphics[width=.45\linewidth]{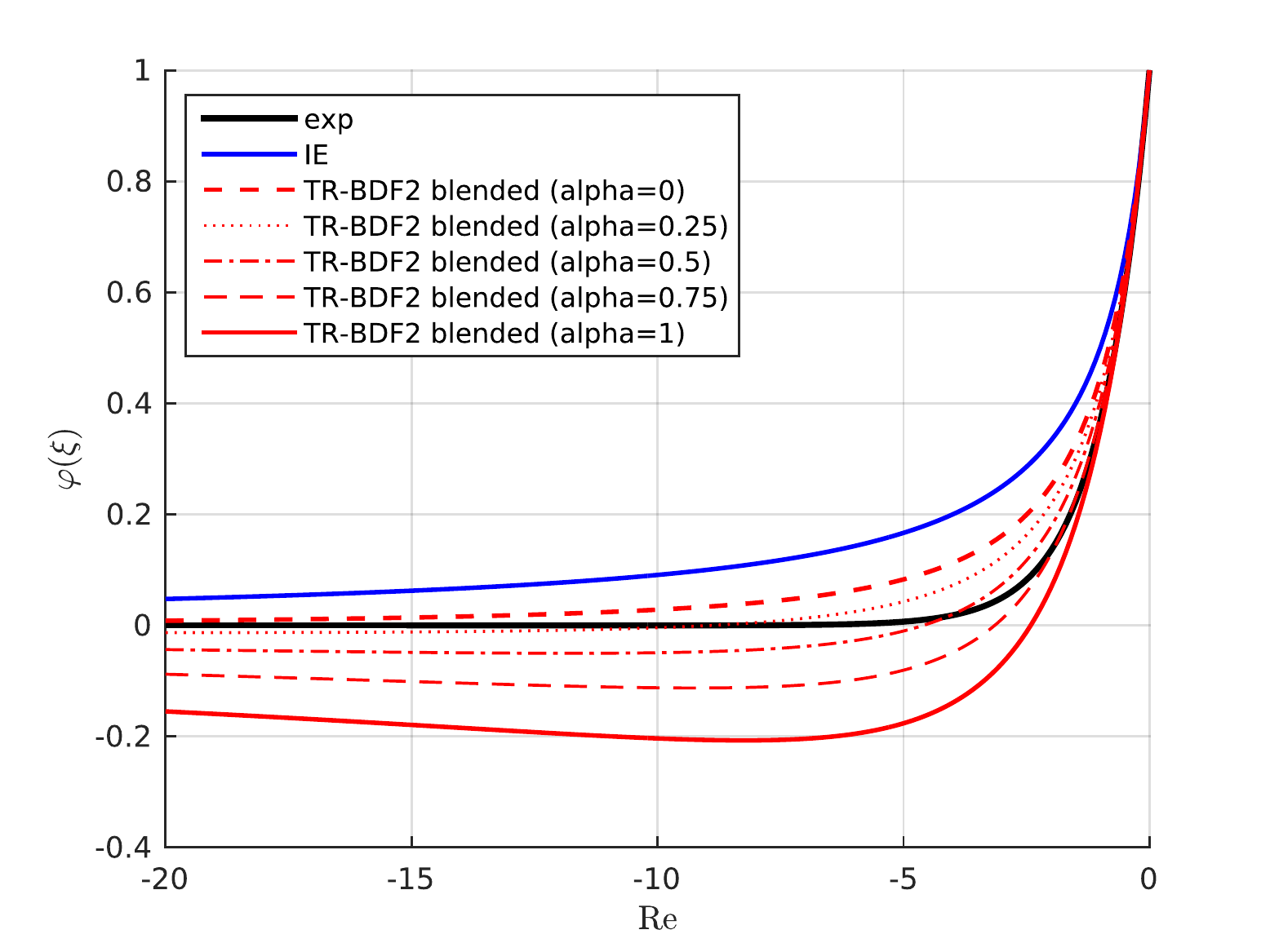}(a)
    \includegraphics[width=.45\linewidth]{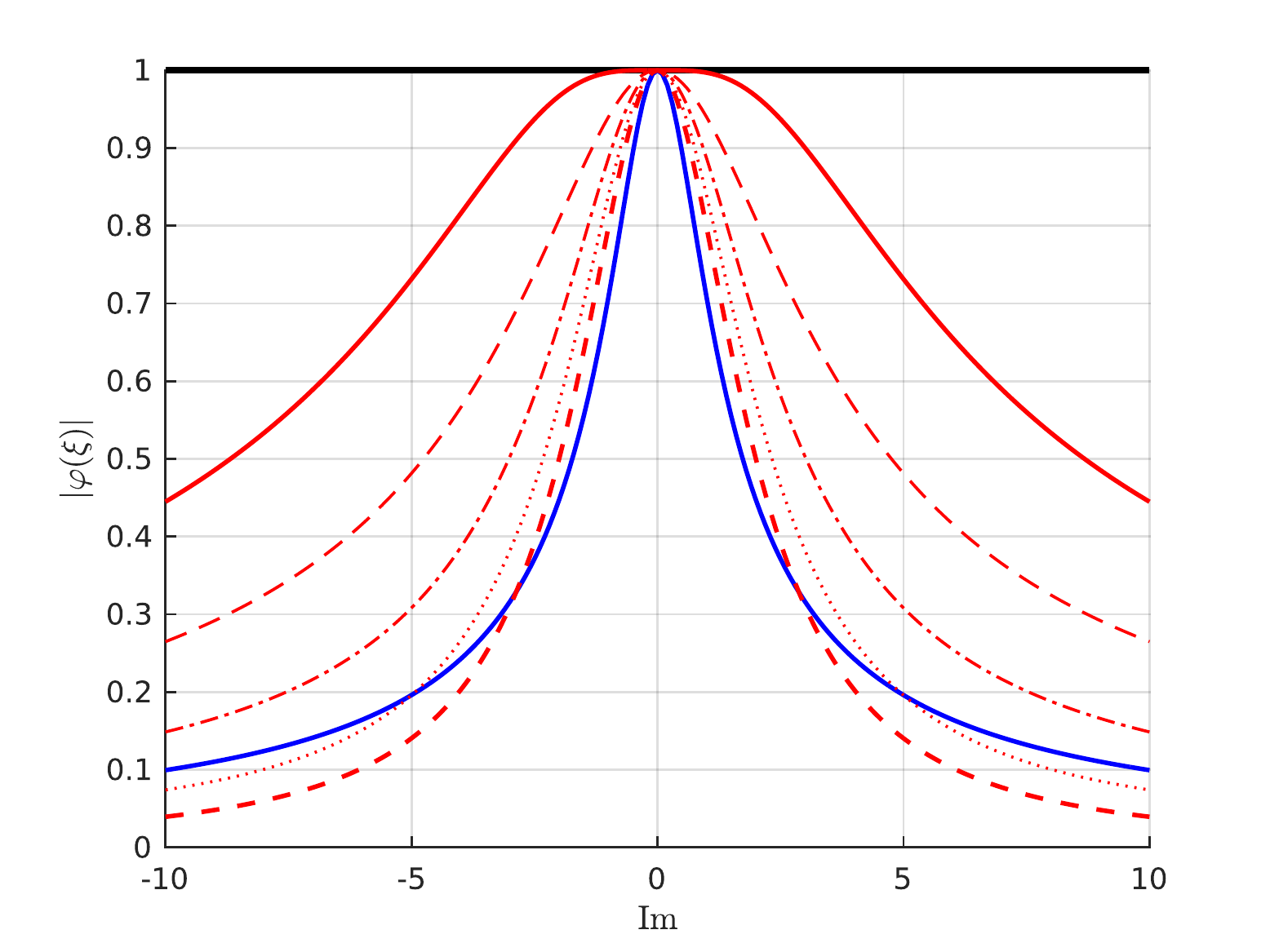}(b)\\
    \includegraphics[width=.45\linewidth]{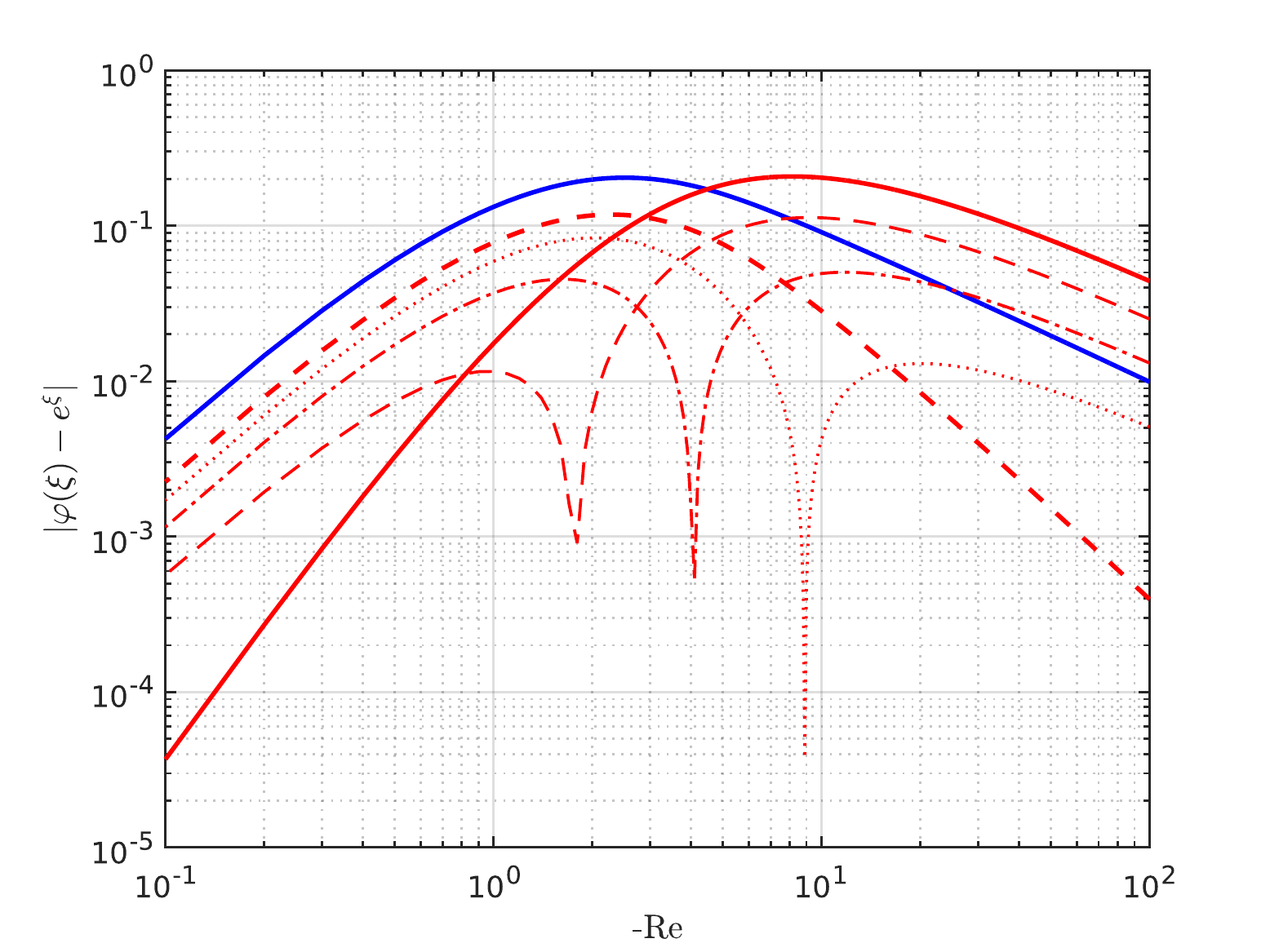}(c)
    \includegraphics[width=.45\linewidth]{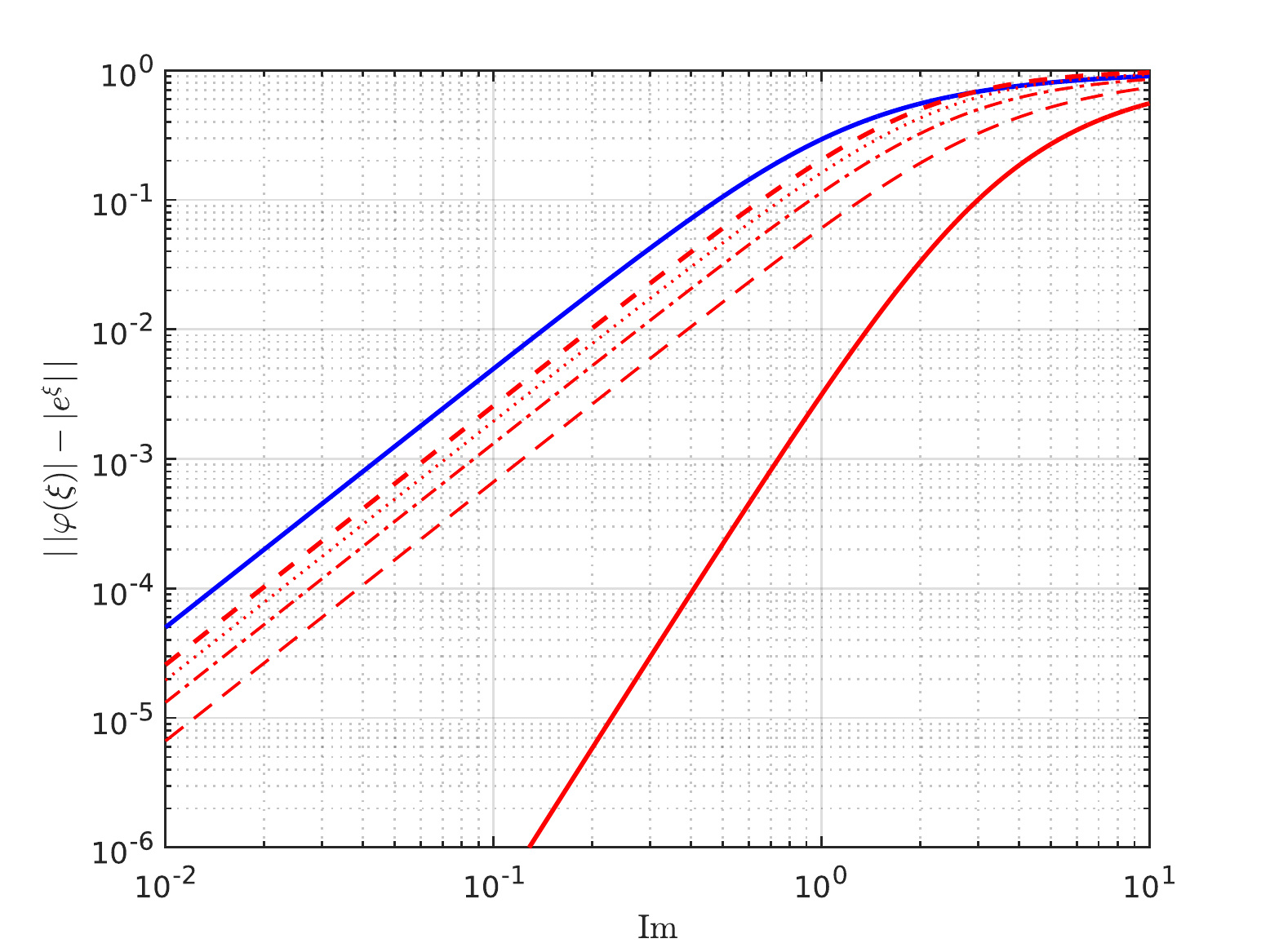}(d)
 \caption{Hybrid TR-BDF2 in comparison with implicit Euler: \emph{a}) stability functions along the negative
   real axis, \emph{b}) modulus of the stability functions along the imaginary axis, \emph{c}) absolute error
   functions along the negative real axis showing the asymptotic and non-asymptotic ranges, \emph{d}) absolute
   error functions along the positive imaginary axis.}
 \label{fig:trbdf2_hybrid}
\end{figure}

The DIRK family corresponding to the hybrid TR-BDF2 method is thus entirely L-stable for every value of the
parameter $\alpha$, as evident from Figure \ref{fig:trbdf2_hybrid}. Additonally for $\alpha \!=\! 1$ the value
$\gamma \!=\! 2-\sqrt{2}$ corresponding to TR-BDF2 maximizes the radius of absolute monotonicity, while for
$\gamma \!=\! 2-\sqrt{2}$ starting from $\alpha \!=\! 1$ (TR-BDF2) the radius of absolute monotonicity is
progressively increased by decreasing $\alpha$, while the order is reduced to $p \!=\! 1$ for $\alpha \!\neq\!
1$, as evident from the absolute error $|\varphi(z) - e^z|$ in the asymptotic range ($Re \!\to\! 0^-$, $Im
\!\to\! 0$) in Figure \ref{fig:trbdf2_hybrid}.  Formally, we have $\lim_{\alpha \to 0} R(A,b) \!=\! \infty$,
thus recovering the unconditional monotonicity of the implicit Euler method.

\begin{figure}[htbp]
    \includegraphics[width=.45\linewidth]{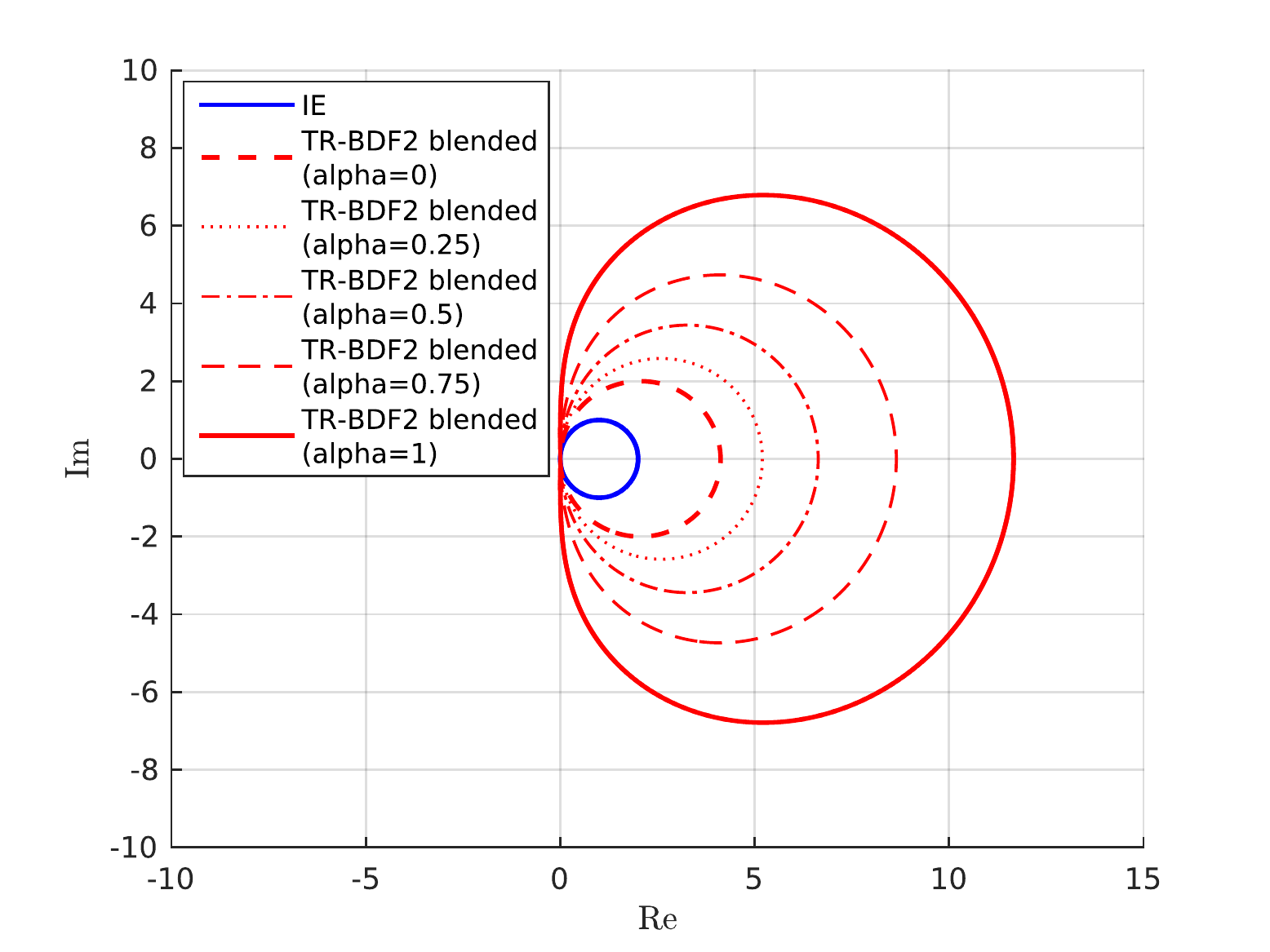}(a)
    \includegraphics[width=.45\linewidth]{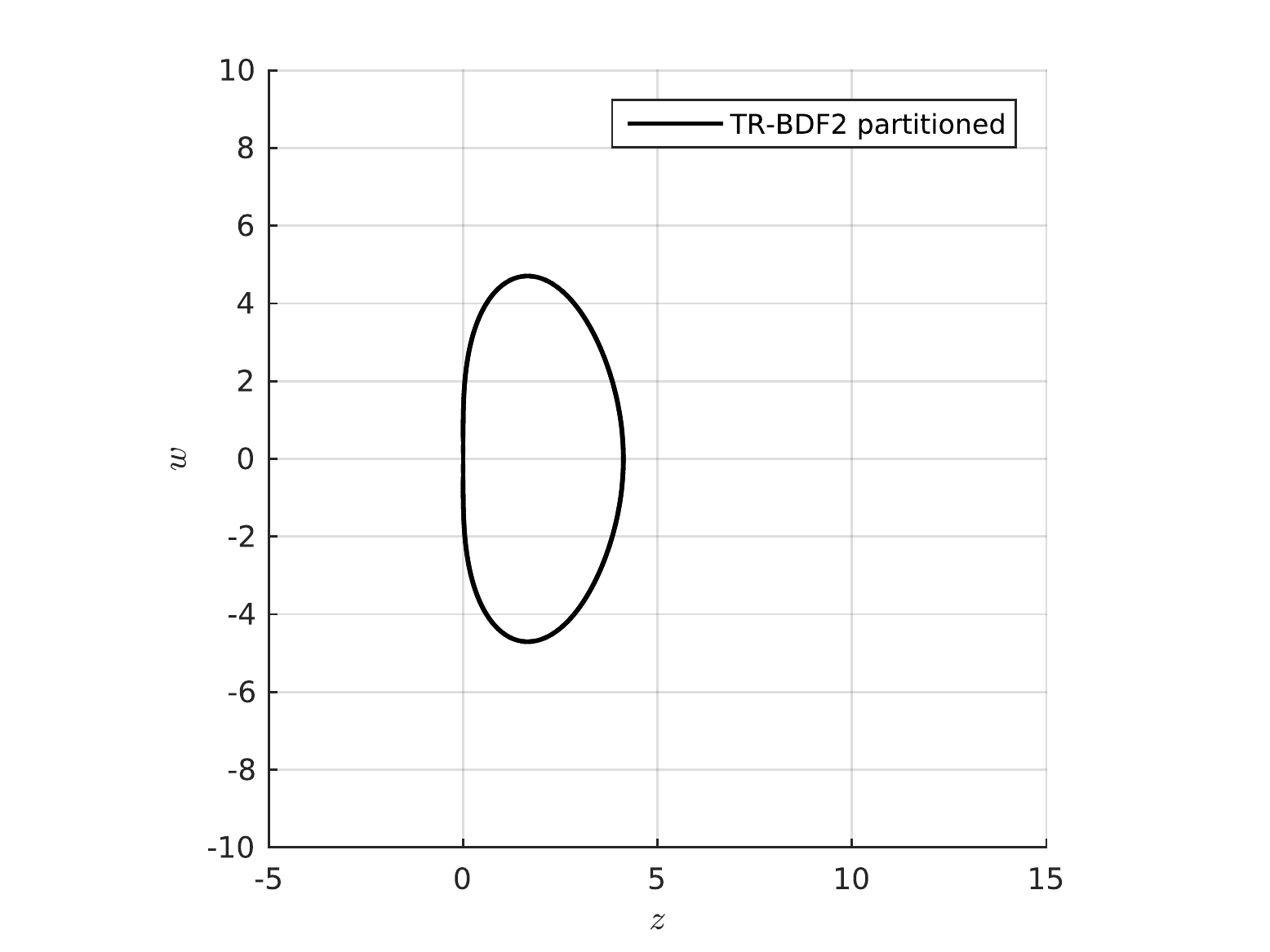}(b)
 \caption{Stability regions (outside portions of the plane): \emph{a}) stability boundaries for the hybrid TR-BDF2
   and the implicit Euler methods, \emph{b}) stability region of the additive Runge-Kutta corresponding to the
   TR-BDF2 partitioned method.}
 \label{fig:trbdf2_hybrid_stability}
\end{figure}

The hybrid TR-BDF2 method can be exploited through different strategies, since by choosing values of the
parameter $0 \leq \alpha \leq 1$ it is possible to produce a continuous blend of the two main schemes
varying the radius of absolute monotonicity accordingly. In our work we adopt a simpler approach and we  
  investigate two alternative modes for enforcing monotonicity under the selected step size.

\subsection{TR-BDF2 blended}

In the \emph{hybrid in time mode} the strategy is to switch from the default $\alpha \!=\!1$ (TR-BDF2) to the
unconditionally monotone mode $\alpha \!=\! 0$ (IE-IE) whenever a suitable sensor detects a violation of some
\emph{global} functional bound on TR-BDF2 solution. After each detected violation by the TR-BDF2 solution,
the time step integration is repeated in IE-IE mode. We call this simple method \textbf{TR-BDF2 blended}, since
it provides automatic adaption of the $\alpha$ value by enforcing unconditional monotonicity only during
critical transients.  To this purpose we introduce the \emph{global} sensor function
\begin{equation}
 \label{eq:global_sensor}
    \sigma = s_g(u^{n+1}) = \begin{cases} 1 & \quad \text{if } \|u^{n+1}\| \leq M \\
                                       0 & \quad \text{otherwise}.\\ \end{cases}
\end{equation}
that is able to determine if the generic functional bound $M$ on $\| \!\cdot\! \|$ is violated by the
numerical solution $u^{n+1}$.  For each time step the algorithm of TR-BDF2 blended with $\gamma \!=\! 2 -
\sqrt 2$ is:
\begin{enumerate}
 \item Set $\alpha \!=\! 1$ and perform the current integration by \eqref{eq:trbdf2_hybrid} to find the tentative
   solution $u^{*}$.
 \item Apply the sensor $\sigma = s_g(u^{*})$.
   \begin{itemize}
    \item If $\sigma \!=\! 1$, set $u^{n+1} = u^{*}$ and go to the next time step.
    \item If $\sigma \!=\! 0$, set $\alpha = \sigma$, repeat the current integration by \eqref{eq:trbdf2_hybrid}
      to find the solution $u^{n+1}$ and go to the next time step.
   \end{itemize}
\end{enumerate}
This basic method is A-stable, L-stable and unconditionally monotone.

\subsection{TR-BDF2 partitioned}

In the \emph{hybrid in solution space mode} the switch from the default TR-BDF2 to the unconditionally
monotone IE-IE mode is applied only for those solution components which are likely to produce violations of a
relevant property under the assigned time step size. In order to detect this behaviour we rely on SSP theory
results by applying Definition \ref{def:discrete_monotonicity} together with the computed value of $R(A,b)
\!\approx\! 2.414$ for TR-BDF2. Since this strategy exploits locality, we introduce the \emph{local} sensor
function
\begin{equation}
 \label{eq:local_sensor}
    \sigma_{i} = s_l(u_{i}^{n+1}) = \begin{cases} 1 & \quad \text{if } \|u_{i}^{n+1}\| \leq M_{i} \\
                                       0 & \quad \text{otherwise}.\\ \end{cases}
\end{equation}
that detects any componentwise violation of the functional bound $M_{i}$ on $\| \!\cdot\! \|$.
In particular for each time step the algorithm is:
\begin{enumerate}
 \item Perform a tentative step of forward Euler
   \begin{equation*}
     u^{*} = u^{n} + h_{EE} \, f(t^n,u^{n})
   \end{equation*}
   using a monotonicity-scaled time step $h_{EE} = h / R(A,b)$.
 \item Apply the sensor $\sigma_i = s_l(u_{i}^{*})$ $(i \!=\! 1,\ldots,m)$ on the tentative solution to construct
   the partitioning matrix $S \!=\! diag\{\sigma_i\}$.
 \item Identifying with $a_{ij}$ and $b_i$ the coefficients corresponding to the tableau for $\alpha=1$ and with
   $\tilde a_{ij}$ and $\tilde b_i$ the coefficients for $\alpha \!=\! 0$ in \eqref{eq:trbdf2_hybrid}, construct
   the automatically \emph{partitioned} RK method
   \begin{subequations}
    \label{eq:trbdf2_partitioned}
     \begin{align}
      g^i &= u^{n} + h \sum\limits_{j=1}^s \,[\,a_{ij} S + \tilde a_{ij} (I-S)\,]\, f(t^n+c_j h,g^j) \\
      u^{n+1} &= u^{n} + h \sum\limits_{i=1}^s \,[\,b_i S + \tilde b_i (I-S)\,]\, f(t^n+c_i h,g^i) 
     \end{align}
   \end{subequations}
   to find the solution $u^{n+1}$ from stage values $g^{i}$ $(i \!=\! 1,\ldots,s)$ with $\gamma \!=\! 2 -
   \sqrt 2$.
\end{enumerate}
%
%
Thus we have effectively transformed the hybrid TR-BDF2 \eqref{eq:trbdf2_hybrid} into a partitioned
Runge-Kutta method, in which the solution space is automatically sorted into monotone and non-monotone
components by performing a tentative forward Euler step on a suitably scaled time step. It may also be
interpreted as an additive Runge-Kutta method (see, e.g., \cite{higueras:2006b}) from which the usual ARK
stability region
$\varphi(z,w) = 1+(i w b^{\intercal} + z \tilde b^{\intercal}) (I -i w A -z \tilde A )^{-1} e$
for the scalar test problem $u^{\prime}(t) \!=\! \lambda u(t) + i \mu u(t)$ with $u(0) \!=\! u^{0}$ is
represented in Figure~\ref{fig:trbdf2_hybrid_stability}. In this case the unconditionally monotone component
$z$ is computed by the IE-IE method while the conditionally monotone component $w$ is integrated using the
original TR-BDF2. We call this second strategy \textbf{TR-BDF2 partitioned}. It introduces a small overhead in
computational time since it always performs an explicit tentative step, even in case of a successful
integration from TR-BDF2.
Clearly, the overall order of accuracy for both strategies will be limited to $p \!=\! 1$ whenever the sensor
functions are activated in the current time step and similarly for the stage order $\tilde
p,$ but both strategies preserve   L-stability and  do not spoil the
linearity of the base methods, they also preserve any linear   invariant and as such they allow atomic
mass conservation (see, e.g., \cite{sandu:2001}).

\section{Potential competitors of TR-BDF2}
 \label{competitors}

In order to compare the properties of TR-BDF2 to those of other similar
 methods in numerical test problems, we introduce here some 
second order methods with similar characteristics. 
We will make our assessment by comparing the performance of TR-BDF2 and its
two hybrid variants from Section~\ref{trbdf2_variants} against the following methods as well as against other
classic methods, such as implicit Euler ($R(A,b)=\infty$) and Crank-Nicolson ($R(A,b)=2$).

\subsection{SSP-optimal SDIRK 2(2)}

The 2-stage SSP-optimal Runge-Kutta method of second order, namely the second order RK with the largest radius
of absolute monotonicity, was discovered in \cite{ferracina:2008}. Later
it was confirmed in \cite{ketcheson:2009}  to be also the optimal 2-stage second order implicit RK through
extensive numerical search. The Butcher tableau of DIRK 2(2) is
%
\begin{center}
\begin{tabular}{c|cc}
          $\frac{1}{4}$ & $\frac{1}{4}$ &  \\
  [1.0ex] $\frac{3}{4}$ & $\frac{1}{2}$ & $\frac{1}{4}$ \\
  [1.0ex] \hline \\
  [-1.5ex]              & $\frac{1}{2}$ & $\frac{1}{2}$
\end{tabular}
\end{center}
As apparent from the tableau, this method consist of two consecutive applications of the implicit midpoint
rule, that is characterized by $R(A,b)=2$, exactly as the Crank-Nicolson method. As a consequence the radius
of absolute monotonicity of the SSP-optimal SDIRK 2(2) is $R(A,b)=4$. Moreover the method inherits all the
properties of the implicit midpoint rule, being A-stable and symplectic, but   it is not L-stable
as TR-BDF2.

\subsection{ROS2}

Among the potential competitors of TR-BDF2 we also include  the Rosenbrock method proposed in
\cite{dekker:1984} and later applied in atmospheric chemistry problems, see e.g. \cite{sandu:2001},
\cite{verwer:1999}. The Butcher tableau is
%
\begin{center}
\begin{tabular}{c|cc|cc}
          $0$ & $0$ &     & $\gamma$ & \\
  [1.0ex] $1$ & $1$ & $0$ & $-2 \gamma$ & $\gamma$ \\
  [1.0ex] \hline \\
  [-1.5ex]    & $\frac{1}{2}$ & $\frac{1}{2}$ & &
\end{tabular}
\end{center}
with $\gamma = 1 + \frac{1}{\sqrt 2}$. In addition to   L-stability, from numerical
experiments in \cite{verwer:1999} it was found that it   also has interesting positivity properties. This
was empirically justified by observing that the stability function is positive along the entire negative real
axis. In \cite{sandu:2001} this observation was extended to the first two derivatives of the stability
function. However, from the framework in Section~\ref{review} it turns out that $R(A,b) \!= 0$, since the third
derivative of the stability function is negative along the negative portion of the real axis. In spite of
this, the fact that up to the second derivative we have positivity for any negative real value can be
interpreted as a sort of \emph{weak} absolute monotonicity. Even though small violations of arbitrary
functionals cannot be a priori excluded at any step size, from numerical experiments in
Section~\ref{experiments} we found an intrinsic resiliency against violations of the TVD property, even from
initial conditions of limited regularity.

\subsection{Modified Patankar Runge-Kutta}

A suitable second order method for stiff chemical problems is the Modified Patankar Runge-Kutta (MPRK)
introduced in \cite{burchard:2003} and later extended to third order accuracy in \cite{formaggia:2011}.  The
MPRK method achieves unconditional positivity through proper weighting of production and destruction terms and
it conserves the quantity $\sum\limits_{j=1}^m u_j, $ which represents the total mass of the system if species
are expressed as mass concentrations.However,  the mass of the atomic species is not conserved. Additionally, due
to the specific form required to the right hand side, it can be used in a PDE framework only by introducing a source
term splitting.

\section{Numerical experiments}
 \label{experiments}

A number of numerical experiments have been carried out, in order to assess the performance of the TR-BDF2
method and its hybrid variants introduced  in Section~\ref{trbdf2_variants} against the other methods described in
Section~\ref{competitors}.  The range of the test cases covers a reactive zero dimensional test problem, here
adapted to the MPRK method, a one dimensional advection problem, an advection diffusion reaction problem for a
mixture of chemical species, as well as two typical nonlinear conservation laws. For PDE tests, we have considered
discontinous initial conditions, in order to show the emergence of critical issues for monotonicity. This
provides the most stringent test, as we experienced from other computations, not reported here, using more
regular initial conditions. In all test cases the nonlinear system associated to implicit RKs was solved using
MATLAB's \texttt{fsolve} with a tolerance of \texttt{tol} = $10^{-10},$ except when otherwise stated. All
computations were performed on a single   Intel$^{\circledR}$ Core\texttrademark i5-2540M (2.60 GHz)
on a laptop with 4 GB RAM running Linux kernel \texttt{3.13.0-24-generic}. We dot not claim that the
error-workload curves reported here are immediately relevant for the selection of step sizes or numerical methods,
 but we consider them as representative of the relative workload expected from the methods assessed.

\subsection{0-D chemical model problem: the Brusselator}

As a first tes,t we use a typical nonlinear chemical kinetics problem. The same issues are shared by all
chemistry modelling problems, which require positivity for each solution component. Consequently, we adopt both
hybridization strategies from Section~\ref{trbdf2_variants} with global and local positivity sensors built
from the floor norm \eqref{eq:functionals_floor_ceil} with $\chi=0$.
In the zero dimensional chemical model problem we measure the maximum error during time integrations from the
$l^{\infty}$-time absolute error norm for the $i$-th species
%
%
\begin{equation}
 \label{eq:absolute_inf_time_norm}
  \| e_i \|_{\infty} = \max_{n=1,\ldots,N_t} \left| u_i^n - \tilde u_i^n \right|
\end{equation}
where $N_t$ is the number of time steps, $u_i^n$ is the solution in the $n$-th time step and $\tilde u_i^n$ is
the reference solution at the same time step obtained with MATLAB \texttt{ode15s} using absolute and relative
error tolerance levels of \texttt{AbsTol}=$10^{-14}$ and \texttt{RelTol}=$10^{-13}$ respectively. In our
numerical results we will refer to the error on species $i \!=\! 1$, which is assumed to be representative of the
problem.

\begin{figure}[htbc]
  \includegraphics[width=0.45\linewidth]{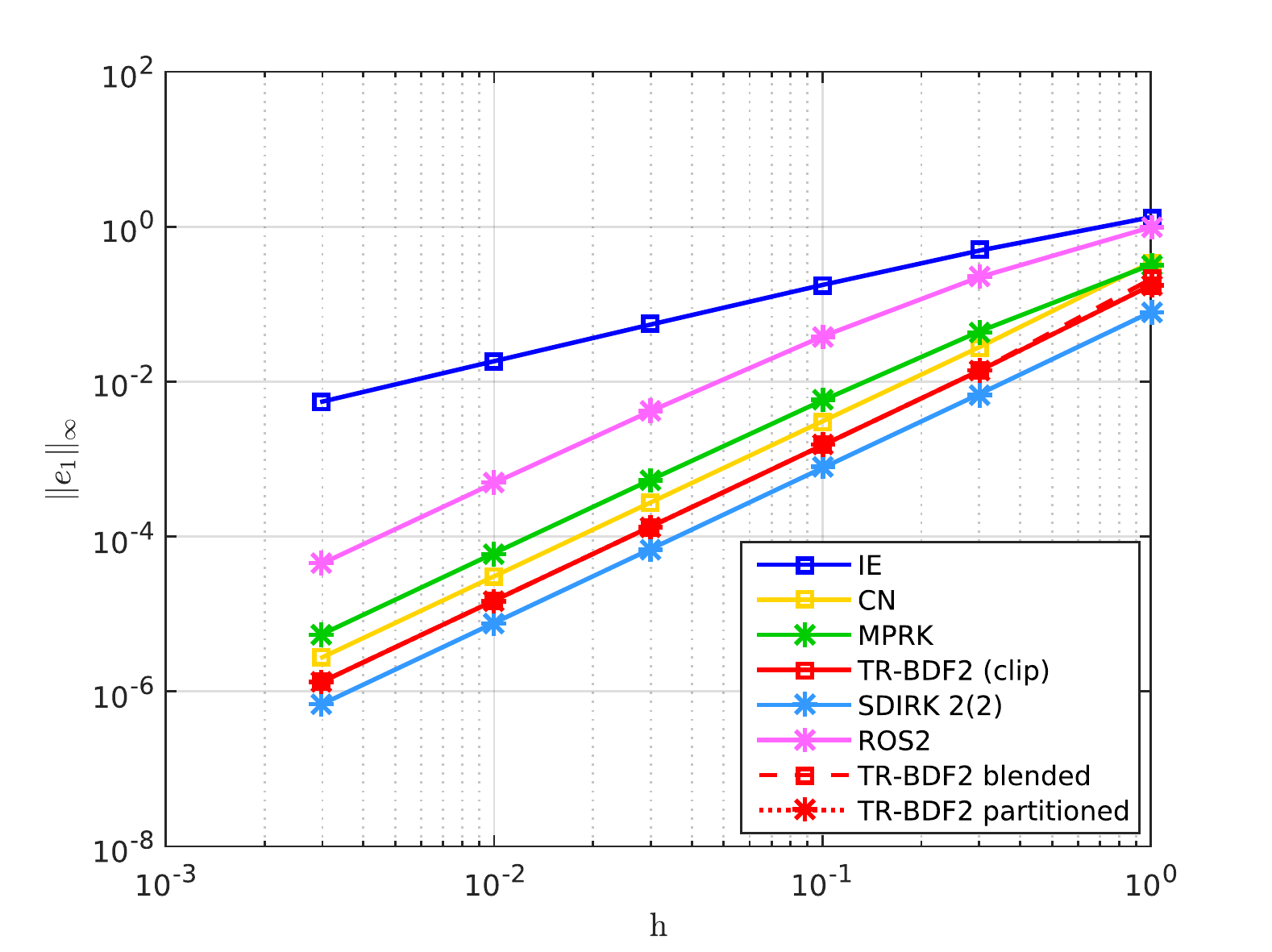}(a)
  \includegraphics[width=0.45\linewidth]{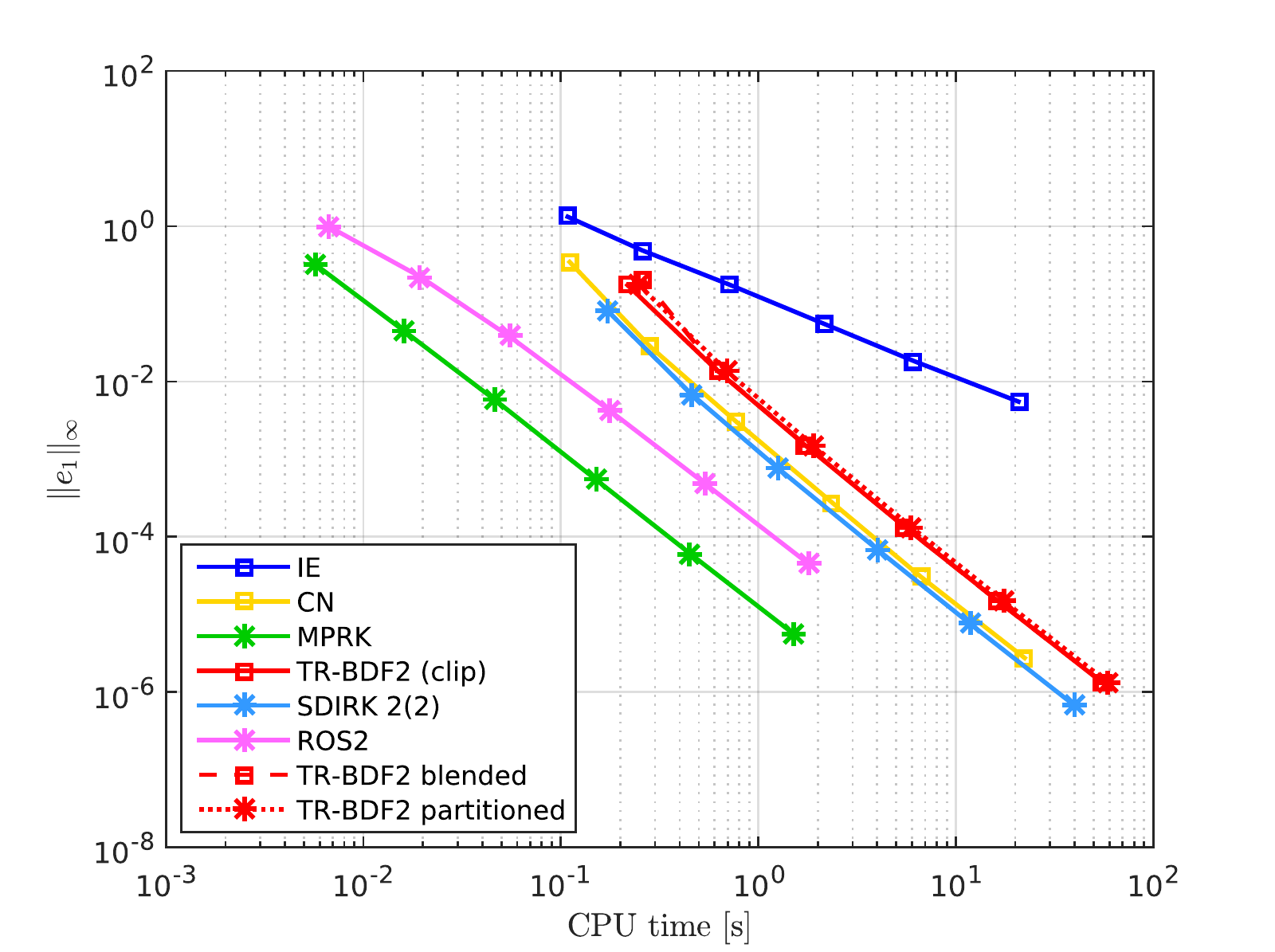}(b)
 \caption{Original Brusselator problem: \emph{a}) error-stepsize curves, \emph{b}) error-workload curves.}
 \label{fig:brusselator_inf}
\end{figure}

We consider thus the original Brusselator system in \cite{lefever:1971}:
\begin{subequations}
 \label{eq:brusselator}
  \begin{align}
   u^{\prime}_1 &= - k_1 u_1 \\
   u^{\prime}_2 &= -k_2 u_2 u_5 \\
   u^{\prime}_3 &= k_2 u_2 u_5 \\
   u^{\prime}_4 &= k_4 u_5 \\
   u^{\prime}_5 &= k_1 u_1 - k_2 u_2 u_5 + k_3 u_5^2u_6 -k_4 u_5 \\
   u^{\prime}_6 &= k_2 u_2 u_5 - k_3 u_5^2 u_6
  \end{align}
\end{subequations}
since this form allows to write the right hand side in a form suitable to MPRK. In particular we follow the
procedure in \cite{formaggia:2011}, which in this case can be applied, since the stoichiometric matrix has
proper rank (\cite{formaggia:2011}, Assumption 2.1). We solve this problem for $0 \leq t \leq 10$ assuming
$k_1 \!=\! k_2 \!=\! k_3 \!=\! k_4 \!=\! k_5 \!=\! k_6 \!=\! 1$ as in the reduced model and starting from the
initial condition $u_1 \!=\! u_2 \!=\! 10$, $u_3 \!=\! u_4 \!=\! 0$ and $u_5 \!=\! u_6 \!=\! 0.1$. The results
in Figure \ref{fig:brusselator_inf} show that Crank-Nicolson, SDIRK 2(2) and TR-BDF2 (clipped to avoid
negative values) are almost equivalent in performance, with a slight advantage for the SSP-optimal SDIRK. MPRK
shows similar accuracy at same step sizes, while it outperforms all the other methods in terms of workload,
being it the only explicit method. ROS2 offers intermediate performance. Blended and partitioned TR-BDF2 are
here equivalent to the clipped version, due to limited size of the integration interval $T$. Similar results,
not shown here, were obtained for the simple geobiochemical problem from \cite{burchard:2003}. Even though the
results from MPRK are promising, in the next tests we are forced to abandon it, since it would require a
source splitting to the advection diffusion reaction problem that is out of our scope.

\subsection{1-D advection problem}

As a first PDE case we consider the  advection equation
\begin{equation}
 \label{eq:test_adv}
  u_t + v u_x = 0, \quad 0 \leq x \leq 1
\end{equation}
which we solve for $0 \!\leq\! t \!\leq\! 1$ and $v \!=\! 1$ using periodic boundary conditions.  We discretize
the interval $(0,1]$ by introducing $N_x$ points $x_i \!=\! i \Delta x$, $i \!=\! 1,\ldots,N_x$ with $\Delta x \!=\! 0.01$
  and we discretize the advective term using a first order upwind scheme to yield a contractive right hand
  side.  We consider the non smooth initial condition
\begin{equation}
 \label{eq:discont_in_cond}
  u(0,x) = \begin{cases} 1 & \quad \text{if } |x-0.5|<0.25\\ 0 & \quad \text{otherwise}.\\ \end{cases}
\end{equation}
For this problem ,the explicit Euler method is stable under the well known condition $\text{Cou} \!=\!
\frac{|a| \Delta t}{\Delta x} \!\leq\! 1$.  In our assessment, we disregard spatial discretization errors,
since we compare the numerical solution obtained with any method to the exact solution of the ODE system
\eqref{eq:ode}, rather than to the exact solution of the original PDE. The exact solution is approximated by
an accurate numerical solution obtained with the MATLAB \texttt{ode45} solver, with absolute and relative error tolerances of
\texttt{AbsTol}=$10^{-14}$ and \texttt{RelTol}=$10^{-13},$ respectively. During numerical tests we measured
%
the $l^{\infty}$-space absolute error norm at final time $t \!=\! T$
\begin{equation}
 \label{eq:absolute_inf_space_norm}
  \| e_i \|_{t=T}^{\infty} = \max_{j=1,\ldots,N_x} \left| u_{i,j}^T - \tilde u_{i,j}^T \right|
\end{equation}
where $u_{i,j}^T$ represents the solution at point $j$ and time $t \!=\! T$ for the single species $i \!=\! 1$
in the one-dimensional advection problem \eqref{eq:test_adv} and $\tilde u_{i,j}^T$ is the reference solution
at the same point in space and time. Furthermore, in order  to assess any violation of the TVD property, we monitored also
the $TV$-space $l^{\infty}$-time seminorm for species $i$
\begin{equation}
 \label{eq:TV_inf_norm}
  \| TV_i \|_{\infty} = \max_{n=1,\ldots,N_t} \sum_{j=1}^{N_x} \left| u_{i,j+1}^n - \tilde u_{i,j}^n \right|.
\end{equation}
%

\begin{figure}[htbc]
  \centering
  \includegraphics[width=.45\linewidth]{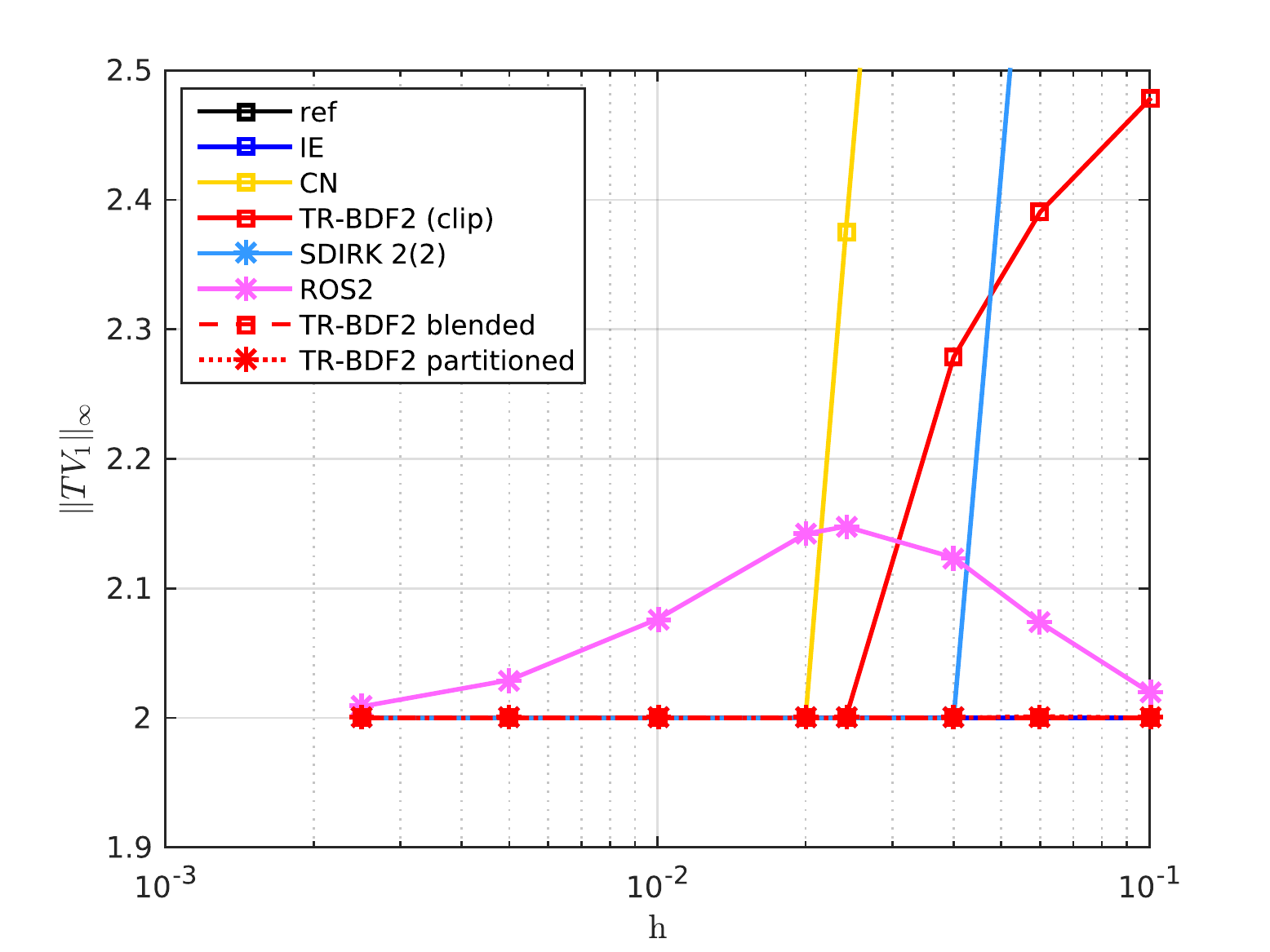}
 \caption{Advection problem with non smooth initial condition: $TV$-space $l_\infty$-time seminorm for the
   single species.}
 \label{fig:adv_upw_1sp_discont_tv}
\end{figure}

\begin{figure}[htbc]
  \includegraphics[width=0.45\linewidth]{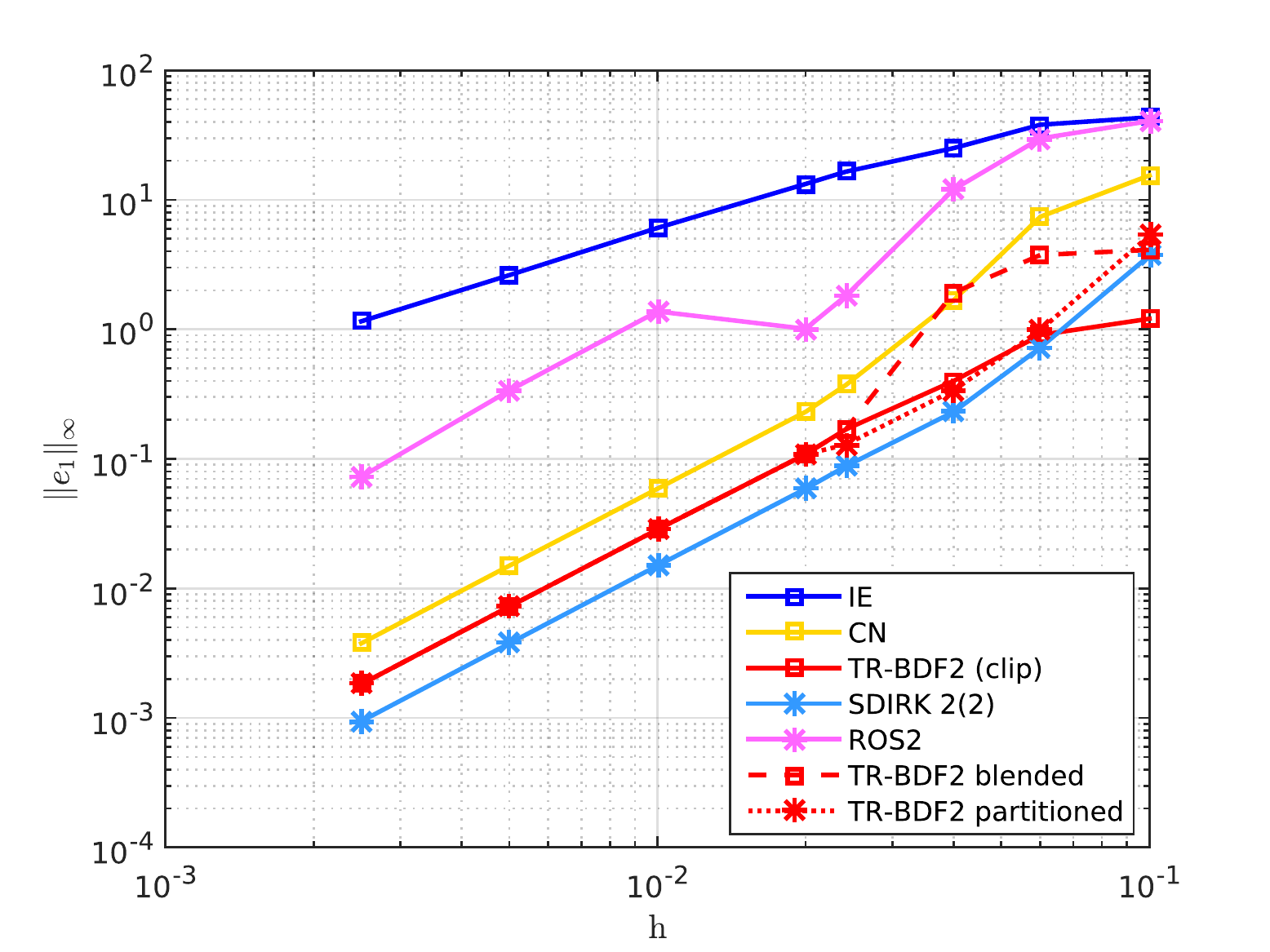}(a)
   \includegraphics[width=0.45\linewidth]{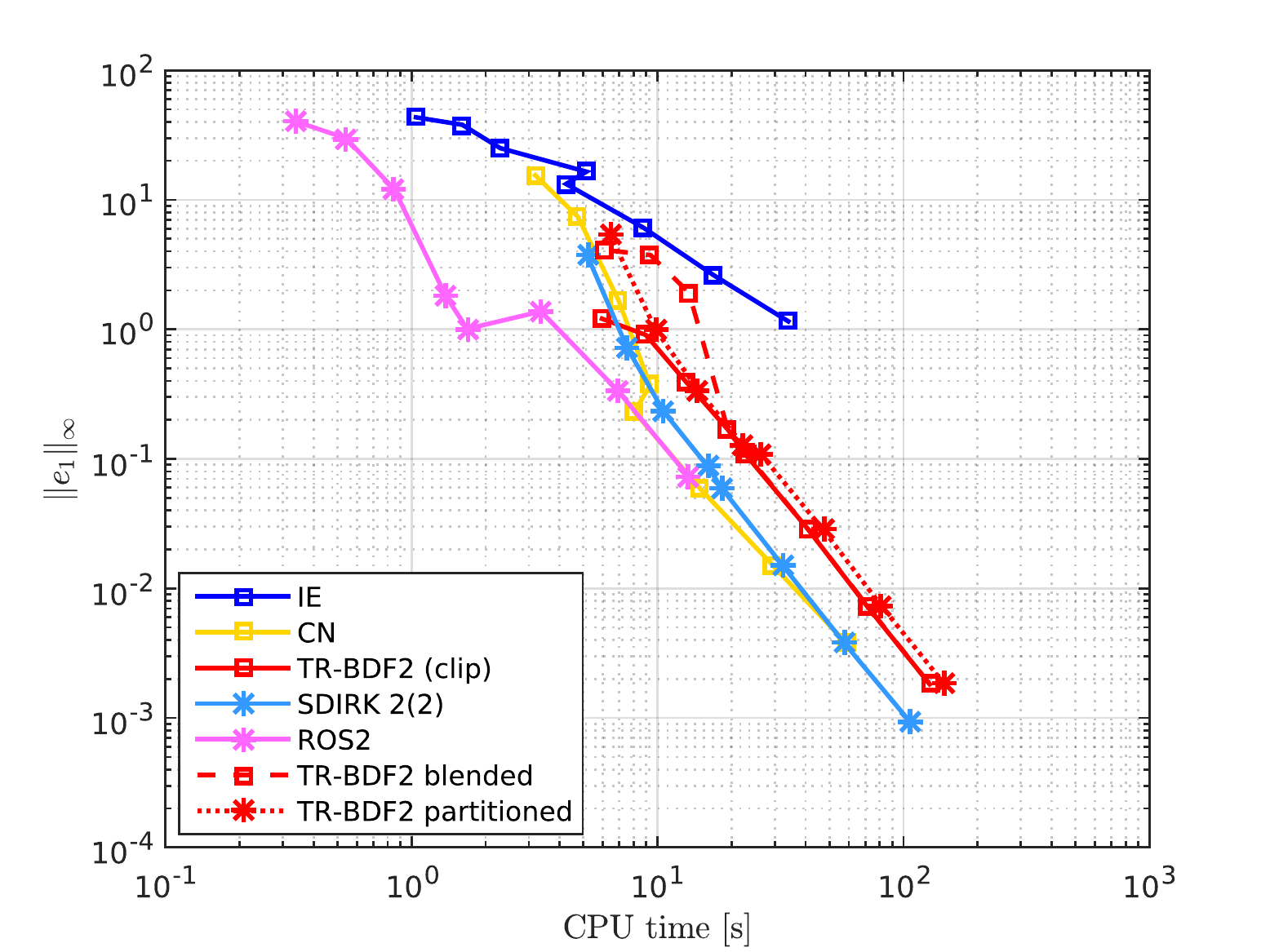}(b)
 \caption{Advection problem with non smooth initial condition: \emph{a}) error-stepsize curves, \emph{b})
   error-workload curves.}
  \label{fig:adv_upw_1sp_discont_inf}
\end{figure}

Here, the TR-BDF2 blended from Section~\ref{trbdf2_variants} relies on a global sensor with the floor
functional with $\chi \!=\! 0$, while TR-BDF2 partitioned uses a local sensor using both floor and ceil
functionals with $\chi \!=\! 0$ and $\psi \!=\! 1$. The discontinuous initial condition
\eqref{eq:discont_in_cond} generates TVD violations for all the conditionally monotone methods (i.e., not for
implicit Euler and TR-BDF2 blended). The critical steps sizes at which violations occur closely
follow the ratio between the radii of absolute monotonicity for the different methods, as from
Sections~\ref{trbdf2} and \ref{competitors}, see Figure \ref{fig:adv_upw_1sp_discont_tv} and Table
\ref{tab:tv_adv_discont_data}. Additionally, ROS2 is never TVD, but the violations are always limited for any
step size and globally it does not show the typical crisis of conditionally monotone methods. This can be
interpreted in light of Section~\ref{competitors}.
TR-BDF2 blended from Section~\ref{trbdf2_variants} never violates TVD property, while its clipped version
shows the usual crisis at large time steps. The possible explanation from Figure
\ref{fig:adv_upw_1sp_discont_inf} is that the critical step size for positivity is smaller than that for TVD
property. As such TR-BDF2 blended starts reverting to the unconditionally monotone mode at relatively small
step sizes, implicitly guaranteeing also the TVD property. On the other hand, TR-BDF2 clipped activates the
clipping procedure for violations of positivity, but it is not able to manage TVD violations. This example
illustrates well the advantage of using SSP theory to guarantee nonlinear solution properties. TR-BDF2
partitioned is always TVD, except for a small violation at the step size $h \!=\! 0.06$.

The accuracy results in Figure~\ref{fig:adv_upw_1sp_discont_inf} maintain the consistent ranking from the zero
dimensional problem. ROS2 now features a non uniform convergence and it achieves higher accuracy for $0.01 \!<\! h
< 0.04$ by sacrificing monotonicity, as evident when comparing Figures \ref{fig:adv_upw_1sp_discont_tv} and
\ref{fig:adv_upw_1sp_discont_inf}.  TR-BDF2 blended shows a degradation in accuracy at larger step sizes, due to
the intervention of the IE-IE mode triggered by the positivity monitor, also evidenced from the order
reduction. The TR-BDF2 partitioned is the only method able to obtain tighter accuracy levels similar to SDIRK
2(2), while additionally mantaining the TVD property with the exception of one step size. When repeating the same
advection test with a smooth initial condition, the results obtained, not shown here, are similar to those in the
chemical model problem and they do not exhibit the order reduction and TVD violations reported above.

\subsection{1-D advection diffusion reaction problem}

As a representative case for chemical transport of reacting species we introduce the advection diffusion
reaction problem of three species
\begin{equation}
 \label{eq:test_adv_diff_react}
  u_t + vu_x = D u_{xx} +  f(u), \quad 0 \leq x \leq 1.
\end{equation}
with periodic boundary conditions. Here $u \!=\! [u_1, u_2, u_3]^T$ and $D \!=\! diag\{d_{ii}\}$ is the
diffusivity matrix. The nonlinear source term $f(u)$ is taken from the simple geobiochemical model in
\cite{burchard:2003}
 \begin{subequations}
 \label{eq:simple_geobio}
  \begin{align}
   f(u)_1 &= - \frac{u_1 u_2}{u_1 +1} \\
    f(u)_2 &= \frac{u_1 u_2}{u_1 +1} - k u_2 \\
   f(u)_3 &= k u_2
  \end{align}
\end{subequations}
representing the interaction among three species identified as nutrients $u_1$, phytoplankton $u_2$ and
detritus $u_3. $ Since the total mass of the system must be conserved, there is an implicit linear invariant
for this problem which is $u_1 + u_2 + u_3$.  We \eqref{eq:test_adv_diff_react} for $0 \leq t \leq 1$ with
$k\!=\!0.3$, an advection velocity of $v=0.1$ and constant diffusivities for each species given by $d_{11}
\!=\! 10^{-3}$, $d_{22} \!=\! 2 \,\text{x} 10^{-3}$ and $d_{33} \!=\! 10^{-4}$.  The grid and the discrete
advection operator are the same of the advection problem. The diffusion term is approximated by central 
finite differencing
and it is naturally contractive.
The initial conditions for the three species are
\begin{subequations}
 \label{eq:adv_diff_react_in_cond}
 \begin{align}
   u_1(0,x) &= \begin{cases} 9.98 \quad &\text{if } |x-0.5|<0.25\\ 0 &\text{otherwise}.\\ \end{cases} \\
   u_2(0,x) &= \begin{cases} 2 \quad &\text{if } |x-0.4|<0.2\\ 0 &\text{otherwise}.\\ \end{cases} \\
   u_3(0,x) &= \begin{cases} 1 \quad &\text{if } |x-0.7|<0.25\\ 0 &\text{otherwise}.\\ \end{cases}
 \end{align}
\end{subequations}
In our assessment, we use a tolerance level \texttt{tol}=$10^{-8}$ for the nonlinear solver in implicit methods
and the reference solution is obtained from MATLAB's \texttt{ode15s} with the same absolute and relative error
tolerances from the advection test. The sensors for the hybrid TR-BDF2 variants are simply built from the
floor functional with $\chi \!=\! 0$.

\begin{figure}[htbc]
   \includegraphics[width=0.45\linewidth]{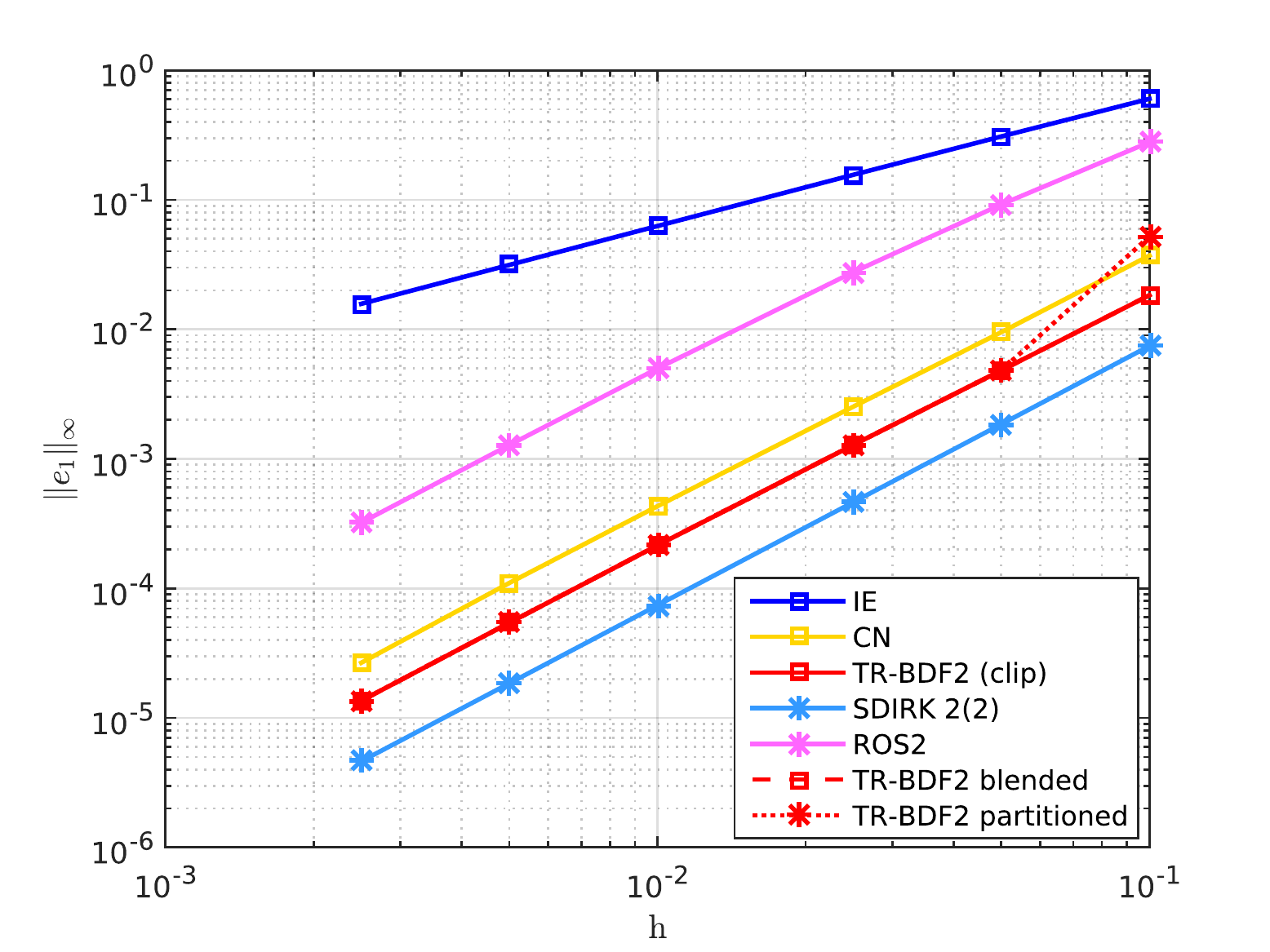}(a)
   \includegraphics[width=0.45\linewidth]{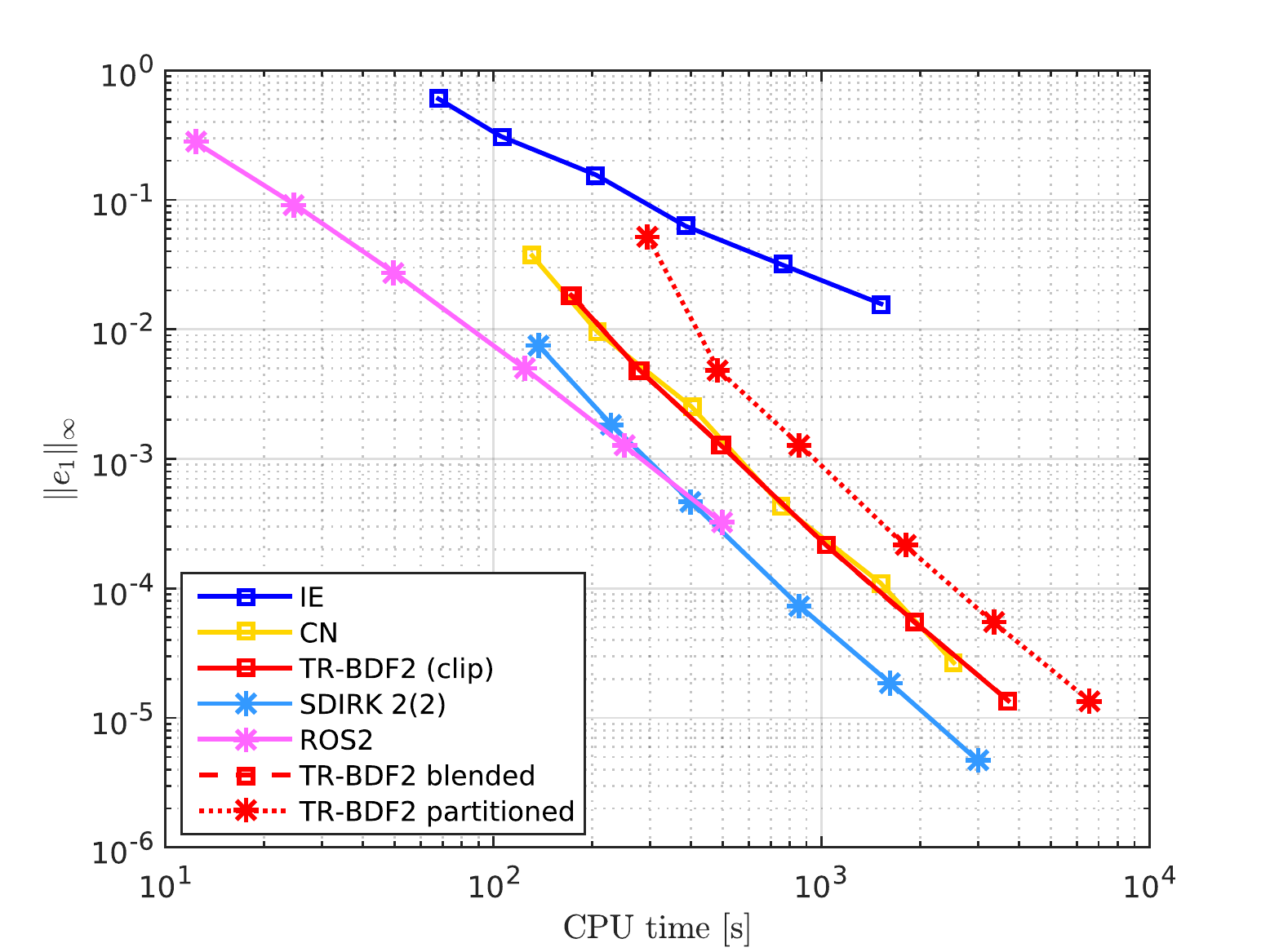}(b)
  \caption{Advection diffusion reaction problem with non smooth initial condition: \emph{a})
   error-stepsize curves, \emph{b}) error-workload curves.}
  \label{fig:adv_diff_react_discont_inf}
\end{figure}

The accuracy curves in Figure \ref{fig:adv_diff_react_discont_inf} do not show the critical features of the
advection test, due to the presence of the additional diffusive term that rapidly smooths out initial
discontinuities. Again, RK methods maintain their relative ranking in terms of accuracy and workload. TR-BDF2 clipped
and the blended variant are almost indistinguishable, while the partitioned version features slightly larger
computational times due to the partitioning step.

\begin{figure}[htbc]
  \includegraphics[width=0.45\linewidth]{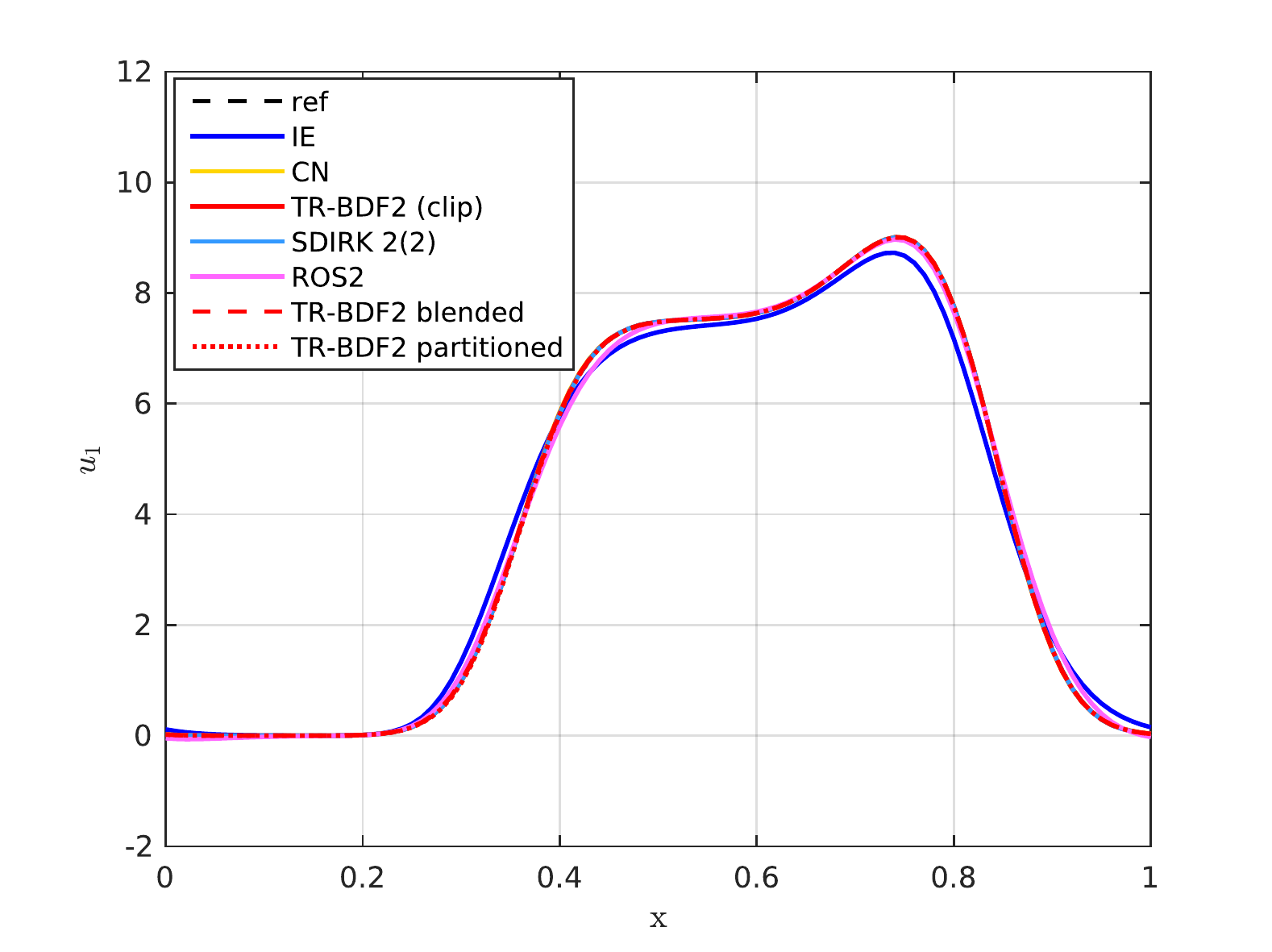}
  \includegraphics[width=0.45\linewidth]{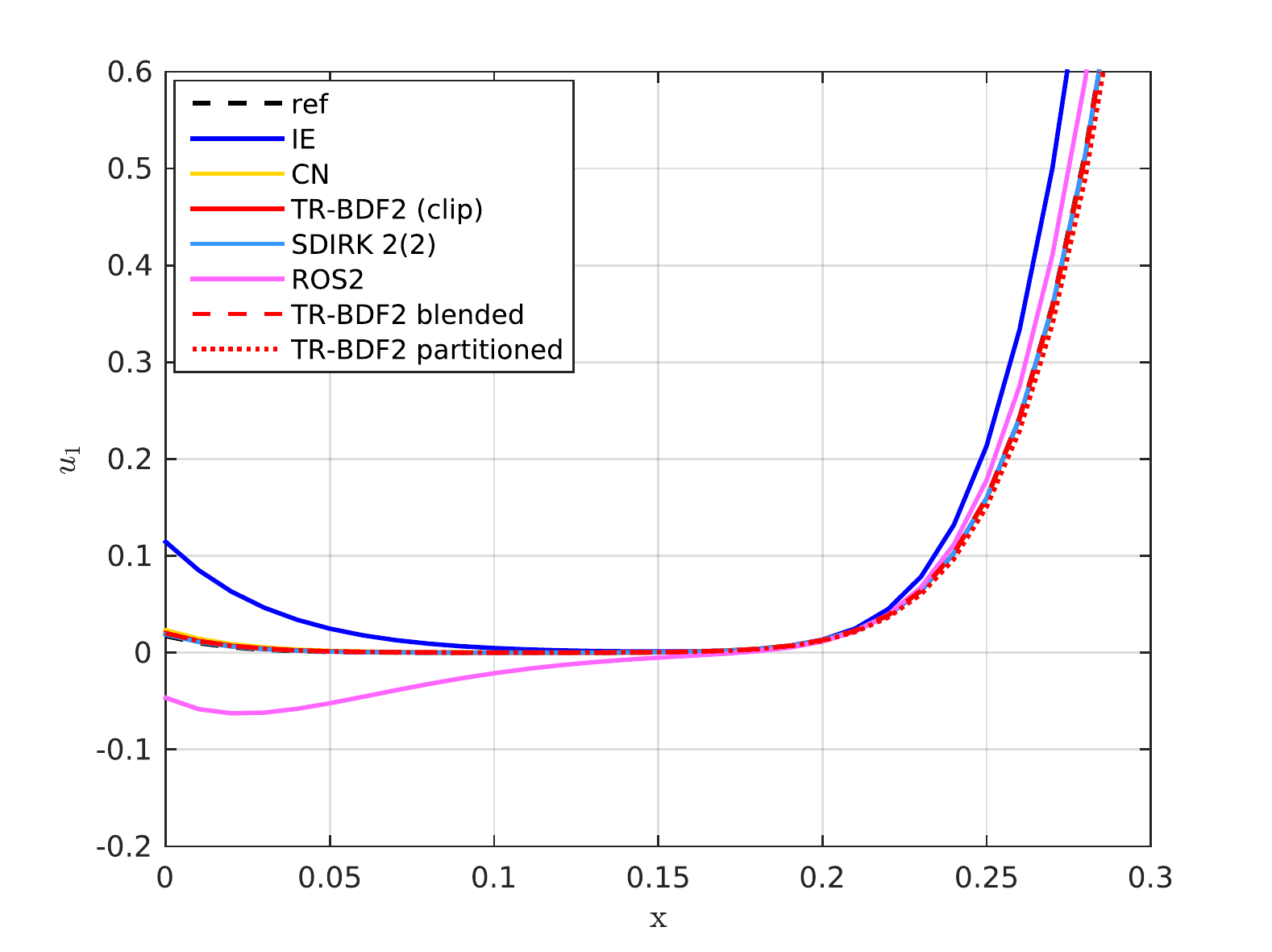}
 \caption{Advection diffusion reaction problem with non smooth initial condition: \emph{a}) solution for
   $u_1(t=T)$ when using $h=0.100$, \emph{b}) close-up in the region of positity violation from ROS2.}
 \label{fig:adv_diff_react_discont_sol}
\end{figure}

\begin{figure}[htbc]
  \centering
  \includegraphics[width=0.45\linewidth]{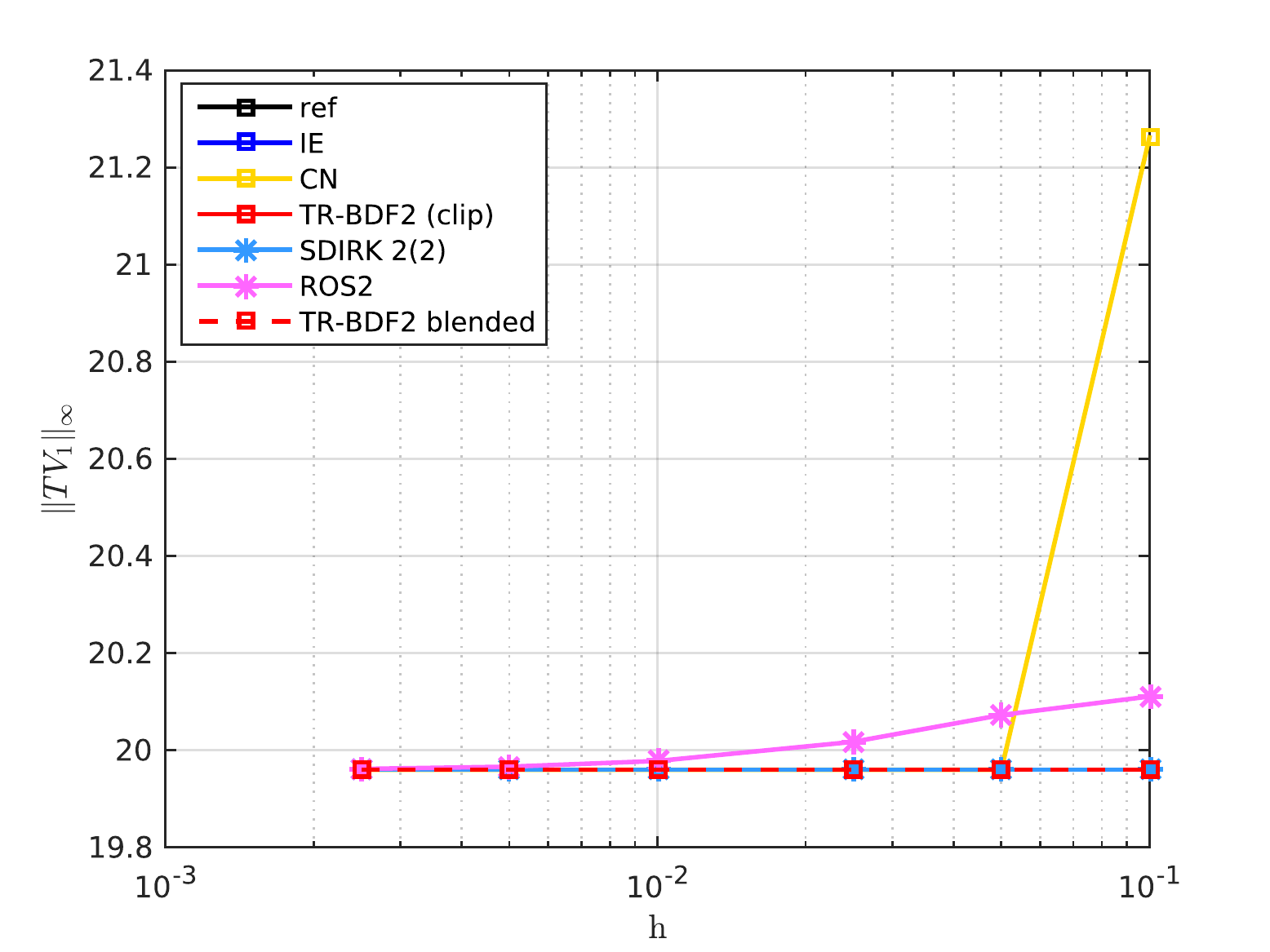}
 \caption{Advection diffusion reaction problem with non smooth initial condition: $TV$-space $l_\infty$-time
   seminorm for the species $u_1$.}
 \label{fig:adv_diff_react_discont_tv}
\end{figure}

A sample solution for this problem is shown in Figure \ref{fig:adv_diff_react_discont_sol}, where a close-up
shows the typical positivity violation on $u_1$ from ROS2. Violations of the TVD property are reported in
Figure \ref{fig:adv_diff_react_discont_tv} where it is evident that all the methods are TVD, with the exceptions
of Crank-Nicolson for $h=0.100$ and ROS2 that is never TVD nor positivity preserving.

\subsection{1-D conservation laws}

We complete our assessment by considering two well known hyperbolic conservation laws, see
e.g. \cite{leveque:2002} for a more detailed discussion. The first is the Burgers equation
\begin{equation}
 \label{eq:test_burgers}
  u_t = -f(u)_x = -\left( \frac{1}{2} u^2 \right)_x, \quad 0 \leq x \leq 1.
\end{equation}
which we solve for $0 \leq t \leq 1$ with the smooth initial condition
\begin{equation}
 \label{eq:smooth_in_cond_burgers}
  u(0,x) = \frac{1}{2} + \frac{1}{4} sin(2 \pi x).
\end{equation}
The second  is the Buckley-Leverett equation
\begin{equation}
 \label{eq:test_burgers}
  u_t = -f(u)_x = -\left( \frac{u^2}{u^2 + \frac{1}{3} (1-u)^2} \right)_x, \quad 0 \leq x \leq 1.
\end{equation}
which we solve for $0 \leq t \leq \frac{1}{8}$ with the discontinuous initial condition
\begin{equation}
 \label{eq:discont_in_cond_buckley_leverett}
  u(0,x) = \begin{cases} \frac{1}{2} & \quad \text{if } x \leq 0.5\\ 0 & \quad
    \text{otherwise}.\\ \end{cases}
\end{equation}
Both equations are here solved with periodic boundary conditions and discretized by 
a high resolution finite volume method using flux
limiters, see e.g.  \cite{hundsdorfer:2003a}, \cite{leveque:2002}.
 More specifically, for the Burgers equation we adopt the \emph{van Leer} limiter
\begin{equation}
 \label{eq:vanleer}
  \Psi(\theta) = \frac{\theta + |\theta|}{1+|\theta|}
\end{equation}
while for the Buckley-Leverett equation we select the \emph{Koren} limiter
\begin{equation}
  \label{eq:koren}
  \Psi(\theta) = \text{max } \left\{ 0 \text{ ; min } \left\{ 2 \text{ ; } \frac{2}{3} + \frac{1}{3} \theta
      \text{ ; } 2 \theta \right\} \right\}.
\end{equation}
Due to the flux limiters, strictly speaking the the Jacobian is not defined. Rather than approximating
it by a finite difference discretization, we again exclude ROS2 from our assessment, since it is known a priori that it performs
poorly on hyperbolic conservation laws.

%

In the numerical tests, we measured the $l^\infty$  norm  in space \eqref{eq:absolute_inf_space_norm}
as well as the $TV$  seminorm in space \eqref{eq:TV_inf_norm}. The reference numerical solution is obtained here
by the MATLAB solver \texttt{ode45} using again \texttt{AbsTol}=$10^{-14}$ and \texttt{RelTol}=$10^{-13}$, while the
implicit stages of the RK methods are solved with a tolerance level of \texttt{tol}=$10^{-10}$. While the TR-BDF2 blended
exploits a global sensor for the TV seminorm, the TR-BDF2 partitioned relies on a local sensor for the floor
and ceil functionals \eqref{eq:functionals_floor_ceil} with $\chi\!=\!0.25$ and $\psi\!=\!0.75$ for the
Burgers equation and $\chi\!=0$ and $\psi\!=\!0.5$ for the Buckley-Leverett equation, respectively.
 These choices follows from the initial conditions \eqref{eq:smooth_in_cond_burgers} and
\eqref{eq:discont_in_cond_buckley_leverett}. Even though this local sensor is not properly a detector of TVD
violations, we use it as an approximate TV sensor due to the known solution dynamics. This is not entirely
correct, as we will see from the tests, but it comes from the difficulty of using a local test for a global
property as TVD.

\begin{figure}[htbc]
  \includegraphics[width=0.45\linewidth]{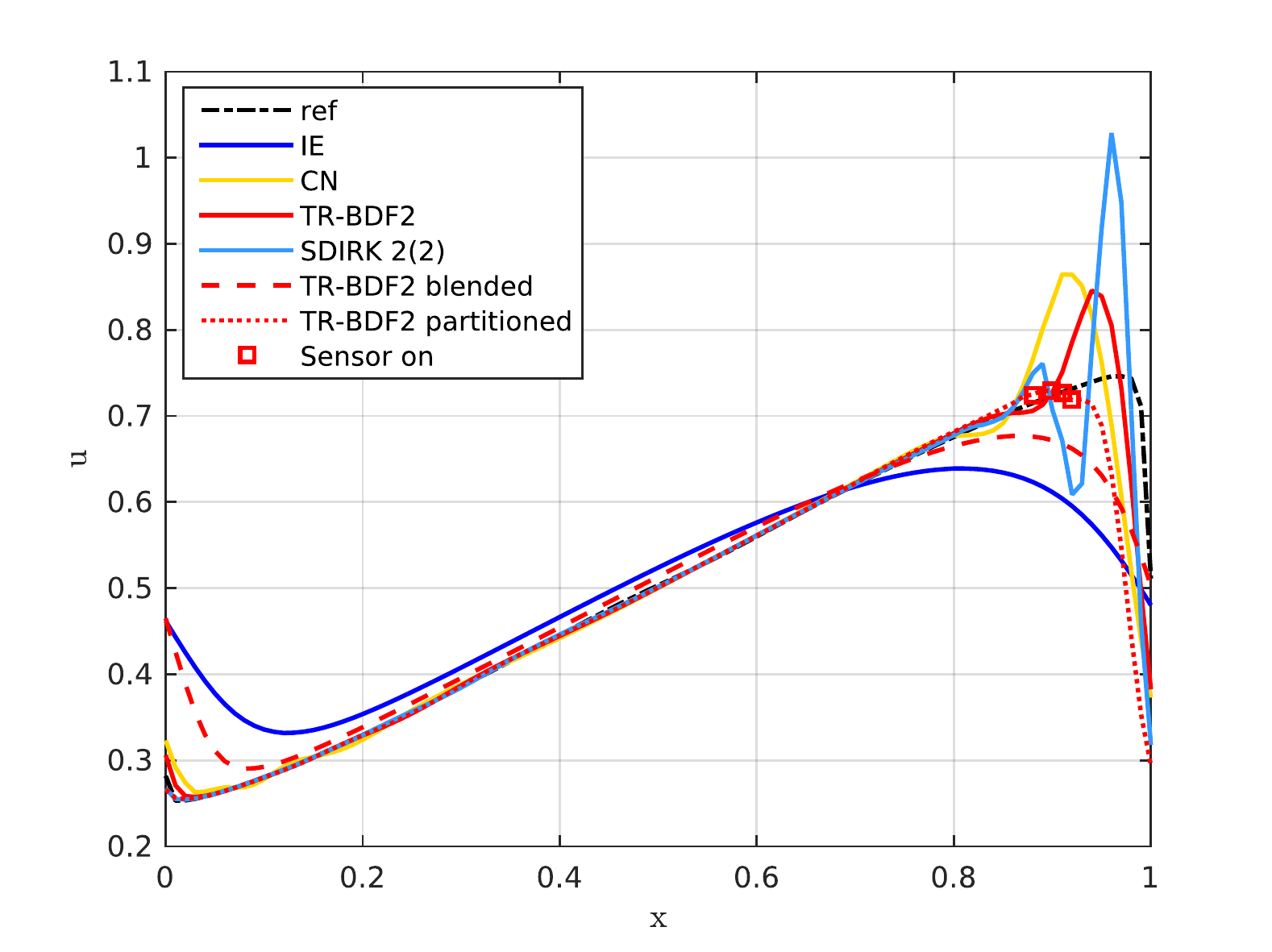}(a)
  \includegraphics[width=0.45\linewidth]{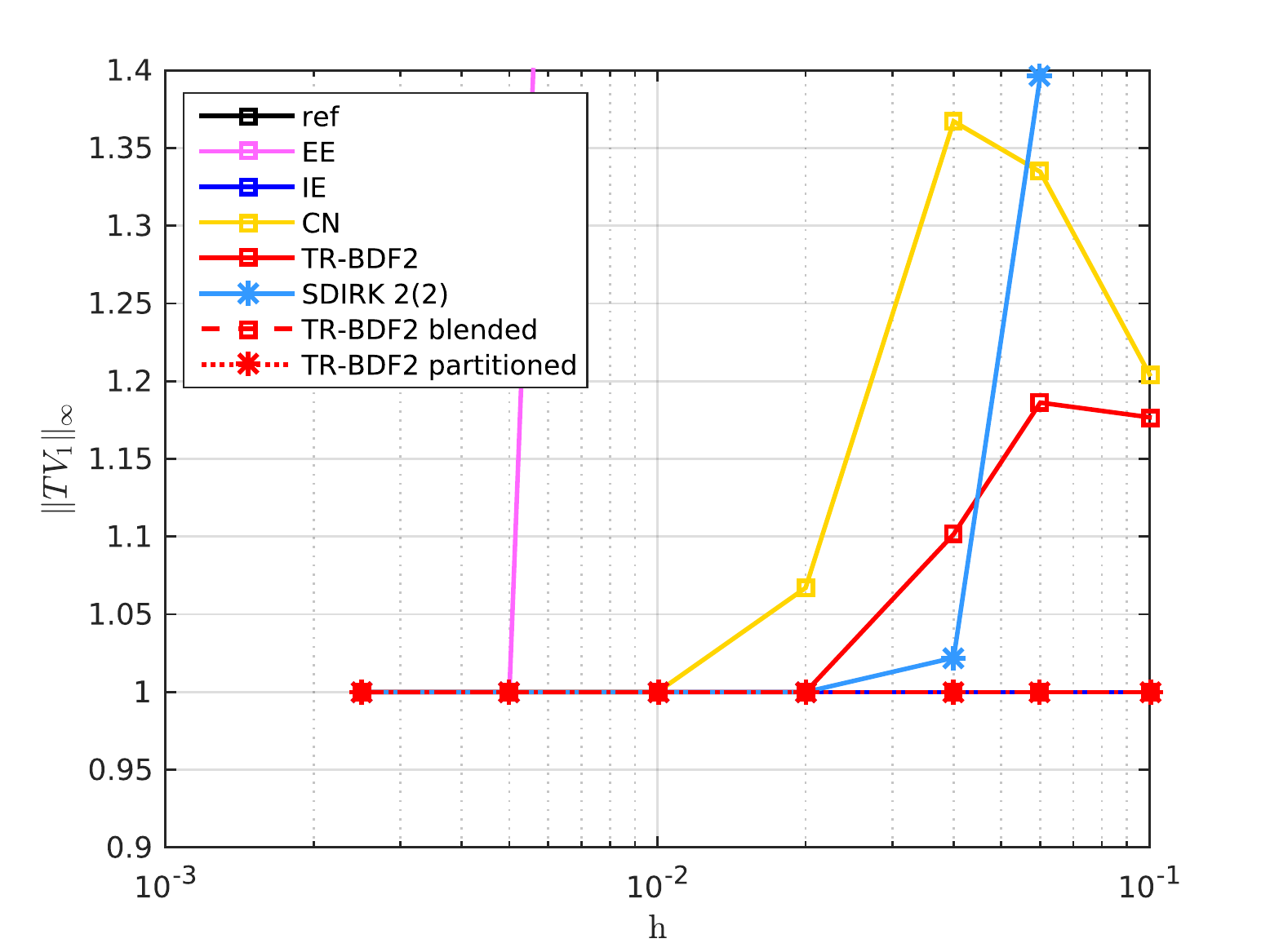}(b)
 \caption{Burgers equation with van Leer limiter and smooth initial condition: \emph{a}) solution at final
   time $t=1$ when using $h=0.1$, \emph{b}) $TV$-space $l_\infty$-time seminorm.}
 \label{fig:burgers_tv}
\end{figure}

\begin{figure}[htbc]
  \includegraphics[width=0.45\linewidth]{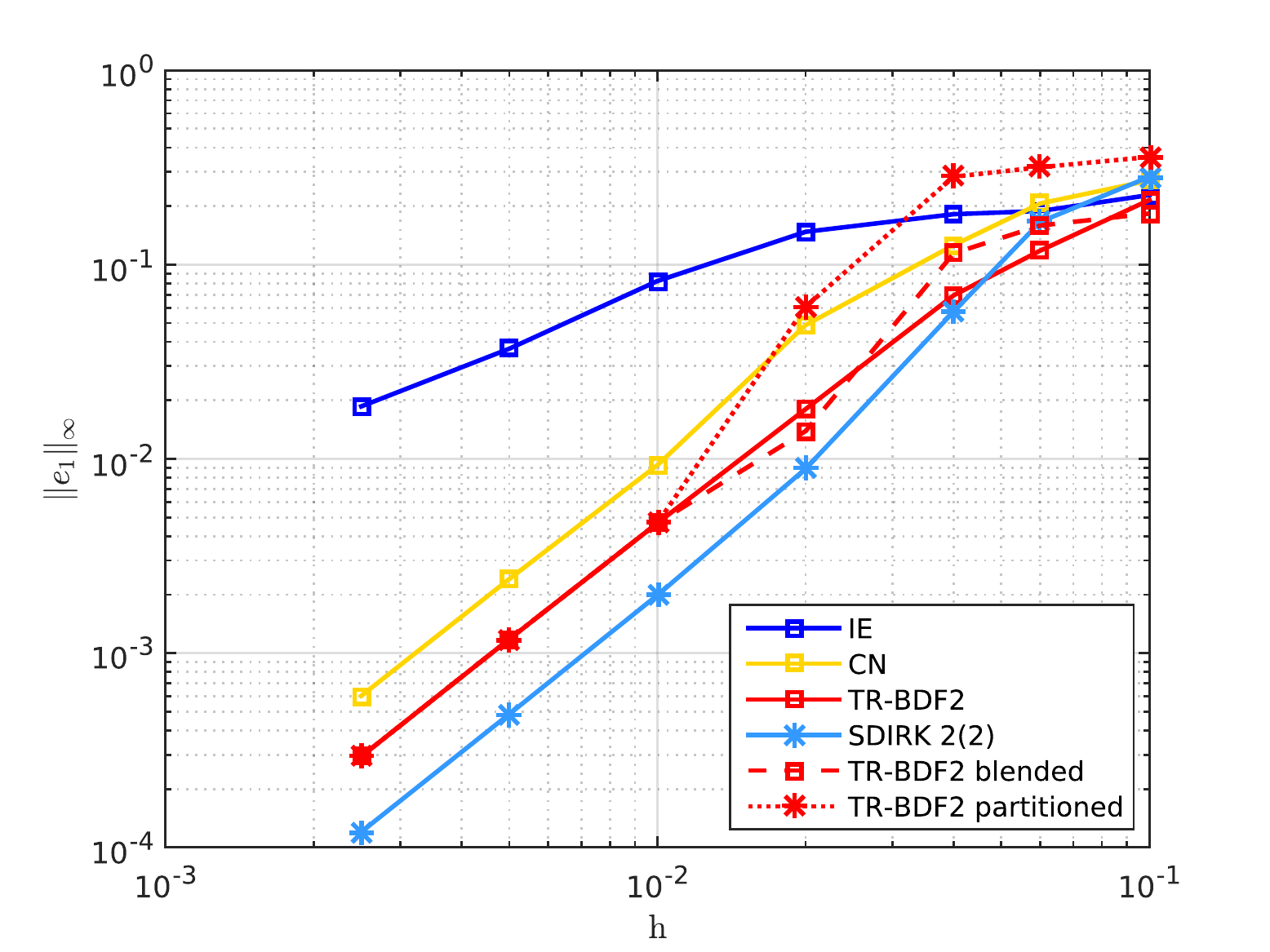}(a)
  \includegraphics[width=0.45\linewidth]{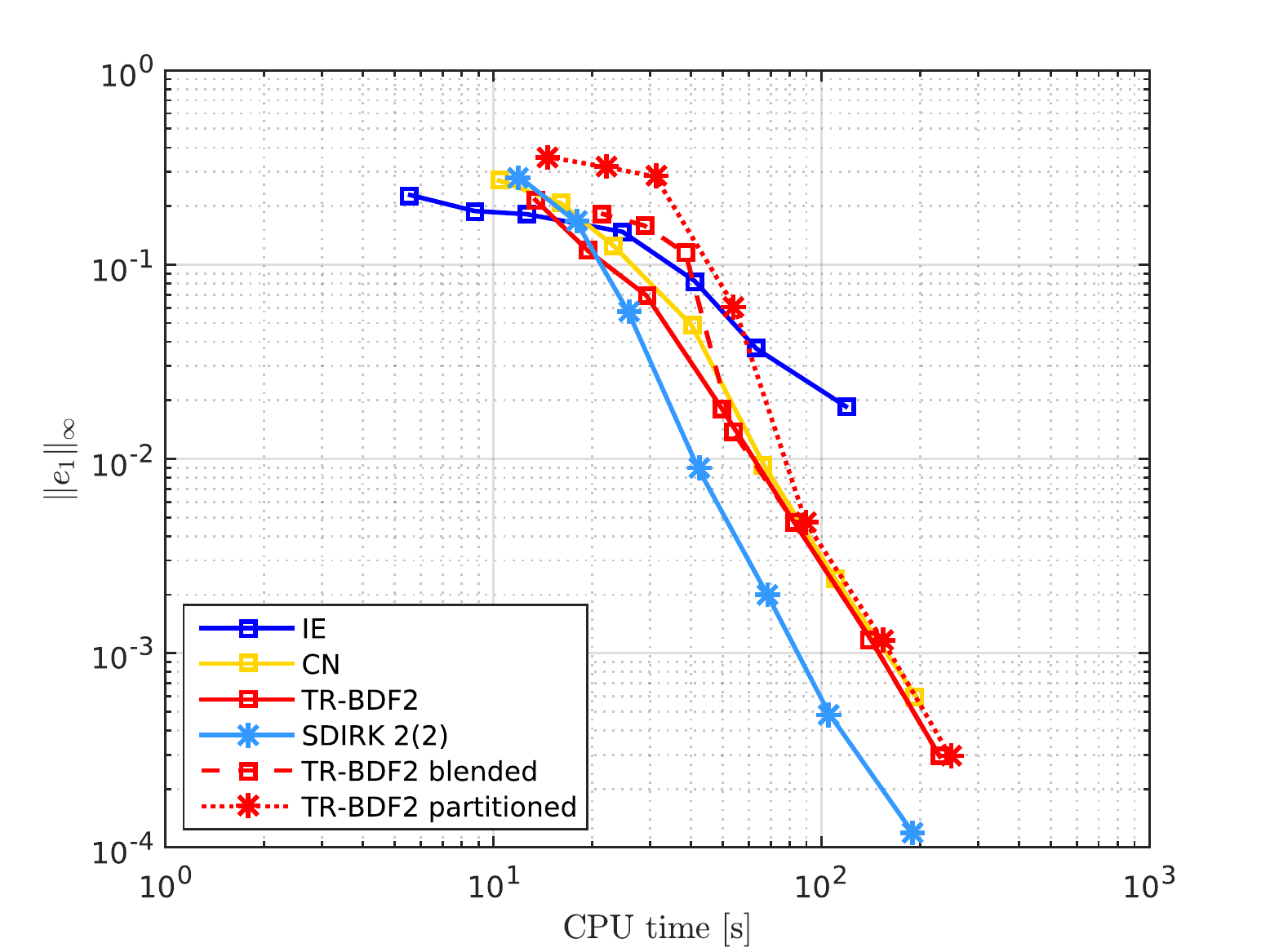}(b)
 \caption{Burgers equation with van Leer limiter and smooth initial condition: \emph{a})
   error-stepsize curves, \emph{b}) error-workload curves.}
 \label{fig:burgers_inf}
\end{figure}

In Figure \ref{fig:burgers_tv} we report the solution at the final time step for the coarsest step size $h
\!=\! 0.1,$ corresponding to about $Cou \!=\! 7.5$. While implicit Euler is able to maintain the TVD property with a
strong smoothing of the developing leading shock, all conditionally monotone methods develop visible
oscillations downstream of this region. In particular, TR-BDF2 (here positive clipping is never activated)
shows minor amplitudes with respect to SDIRK 2(2) and Crank-Nicolson. TR-BDF2 blended obtains a smoothed
solution after several integrations in IE-IE mode. TR-BDF2 partitioned is qualitatively very close to the
reference solution with the best approximation for the shock amplitude, but it features also a reduction in
the propagation speed, probably due to the smoothing of the imaginary parts for the points integrated in IE-IE
mode (see Figure \ref{fig:trbdf2_hybrid}). Significantly, both variants  always  remain TVD, even at $Cou \!=\!
7.5$. The other methods show TVD violations with the usual critical step size progression, see Figure
\ref{fig:burgers_tv}.
The accuracy results for the Burgers equation are reported in Figure \ref{fig:burgers_inf}. Interesting
behaviour appears at coarse stepsizes where all the methods collapse about at the same accuracy of implicit
Euler. TR-BDF2 blended realizes a smooth adaption from the monotone implicit Euler accuracy to the TR-BDF2
asymptotic curve. It preserves always monotonicity, while Crank-Nicolson, SDIRK 2(2) and TR-BDF2 clipped
violate it, as evident from Figure \ref{fig:burgers_tv}. The error curves from TR-BDF2 partitioned follows tha
same behaviour even though the $l^\infty$ error norm is penalized from the behaviour at the leading shock.

\begin{figure}[htbc]
  \includegraphics[width=0.45\linewidth]{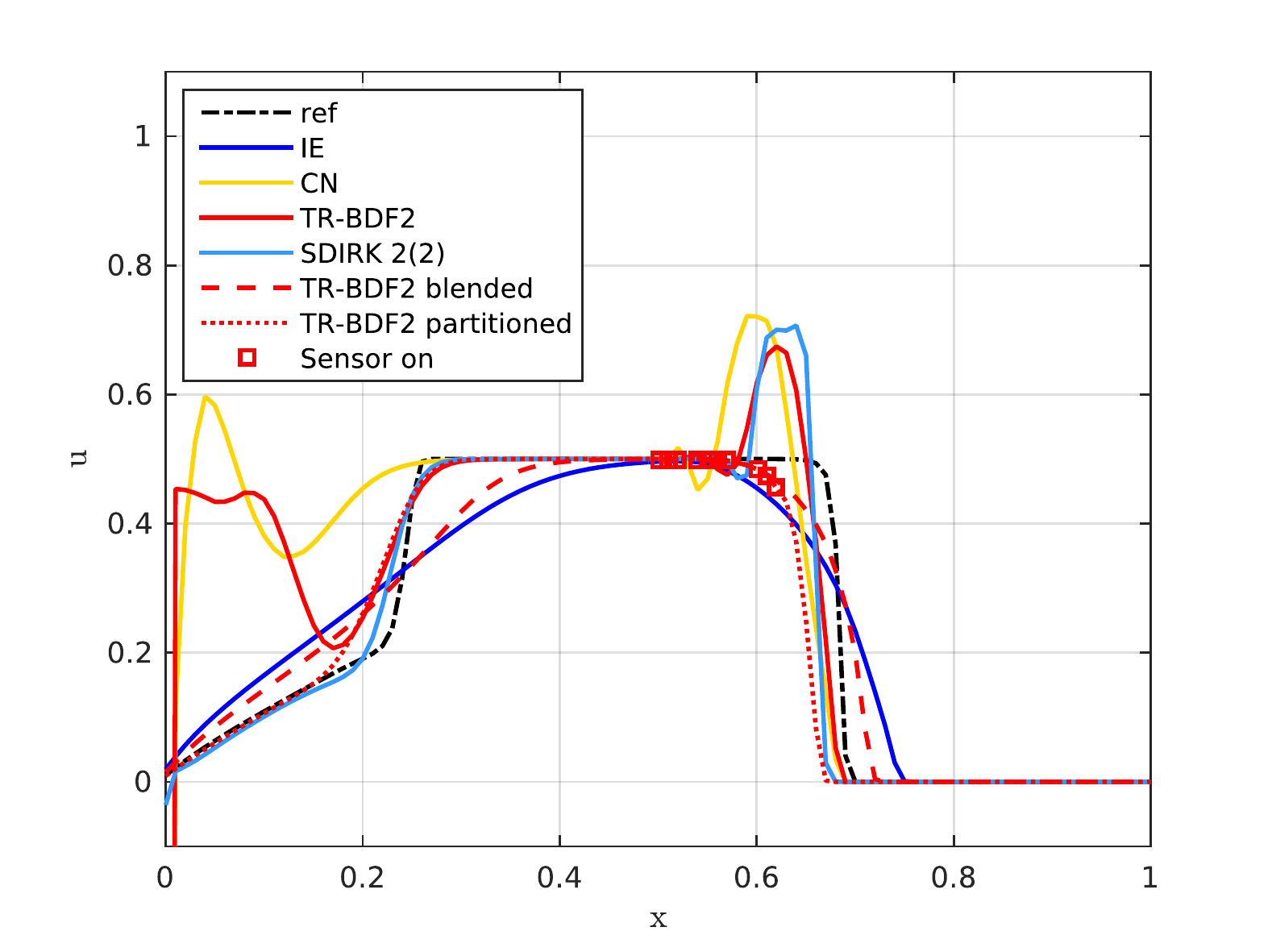}(a)
  \includegraphics[width=0.45\linewidth]{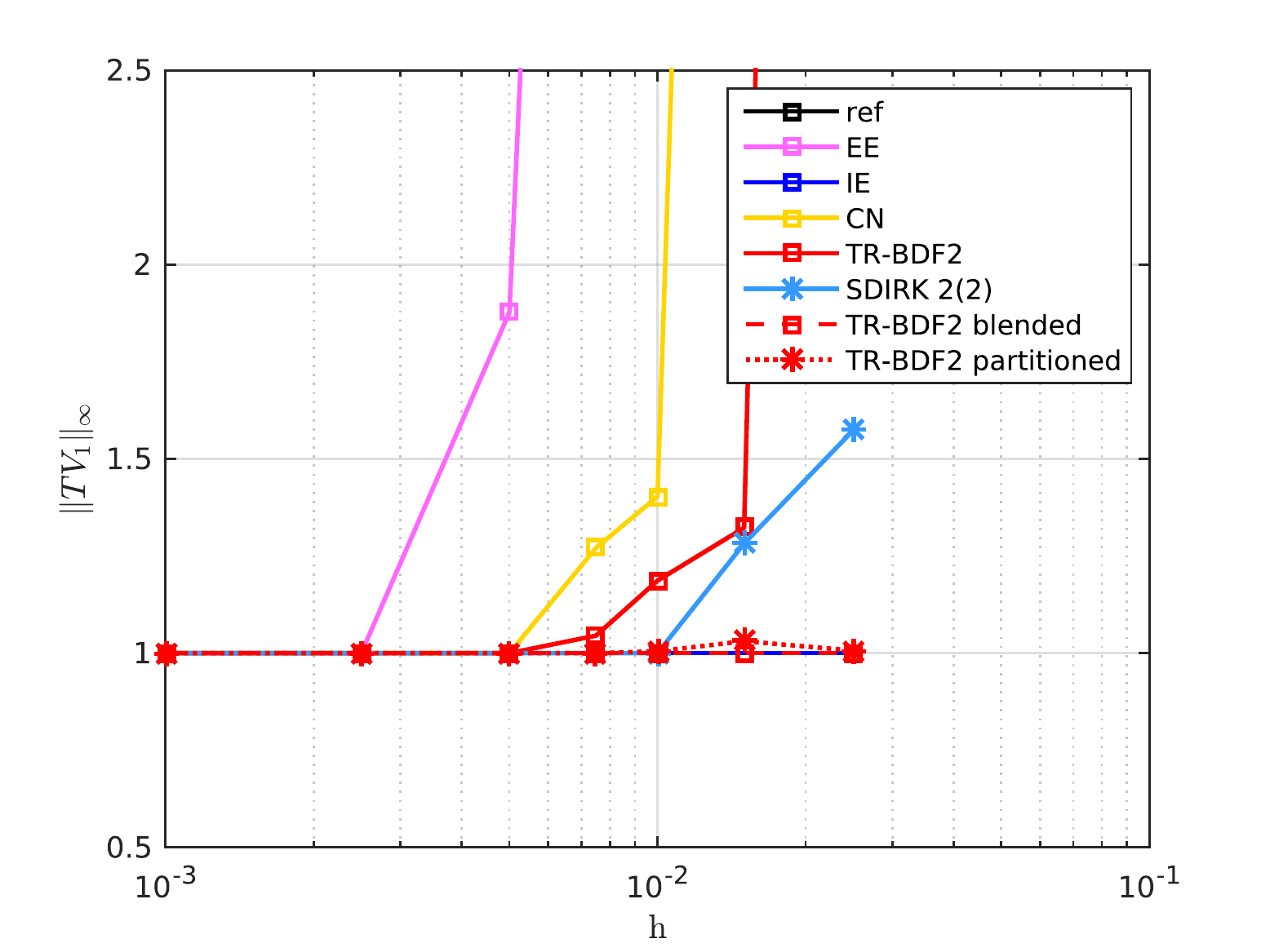}(b)
 \caption{Buckley-Leverett equation with Koren limiter and non smooth initial condition: \emph{a}) solution at
   final time $t=0.125$ when using $h=0.025$, \emph{b}) $TV$-space $l^\infty$-time seminorm.}
 \label{fig:buckley_leverett_tv}
\end{figure}

The Buckley-Leverett test provides a more stringest test due to its non-convex flux function. The solution at
the final time step for the step size $h \!=\! 0.025$ is shown in Figure \ref{fig:buckley_leverett_tv}. Again
all conditionally monotone methods develop oscillation on the trailing shock, while only SDIRK 2 (2) develops
a stable rarefaction wave, where Crank-Nicolson and TR-BDF2 (again without clipping) develop additional
waves. Implicit Euler and TR-BDF2 blended feature a strongly smoothed behaviour, while TR-BDF2 partitioned
remains free of oscillations due to the activation of the local sensor that allows to develop the shock and
rarefaction waves. Anyway it shows a reduction in the shock speed, similarly as in the Burgers test, and some
lmited violations of the TVD property at the largest step sizes, see Figure \ref{fig:buckley_leverett_tv}.

\begin{figure}[htbc]
  \includegraphics[width=0.45\linewidth]{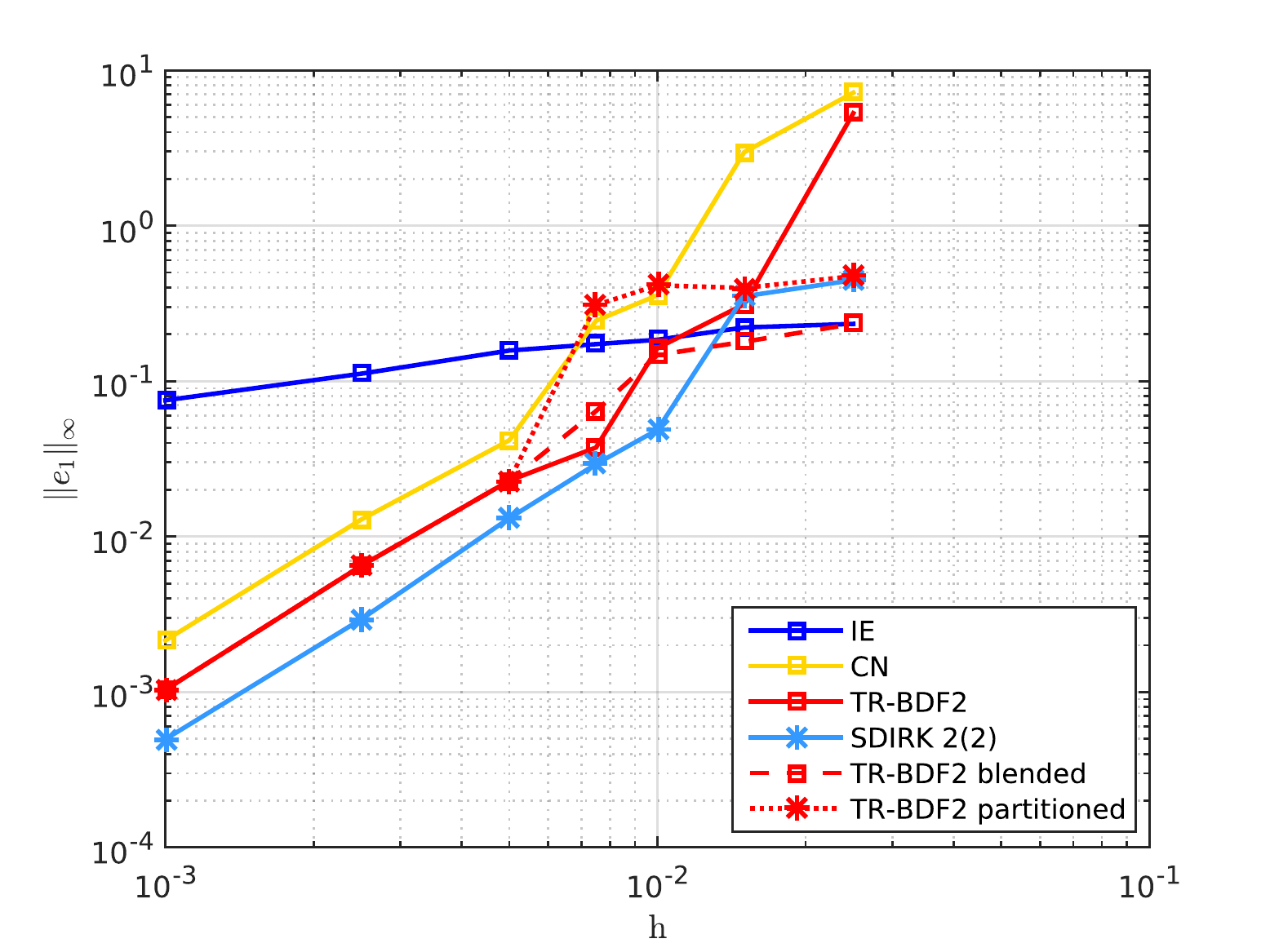}(a)
  \includegraphics[width=0.45\linewidth]{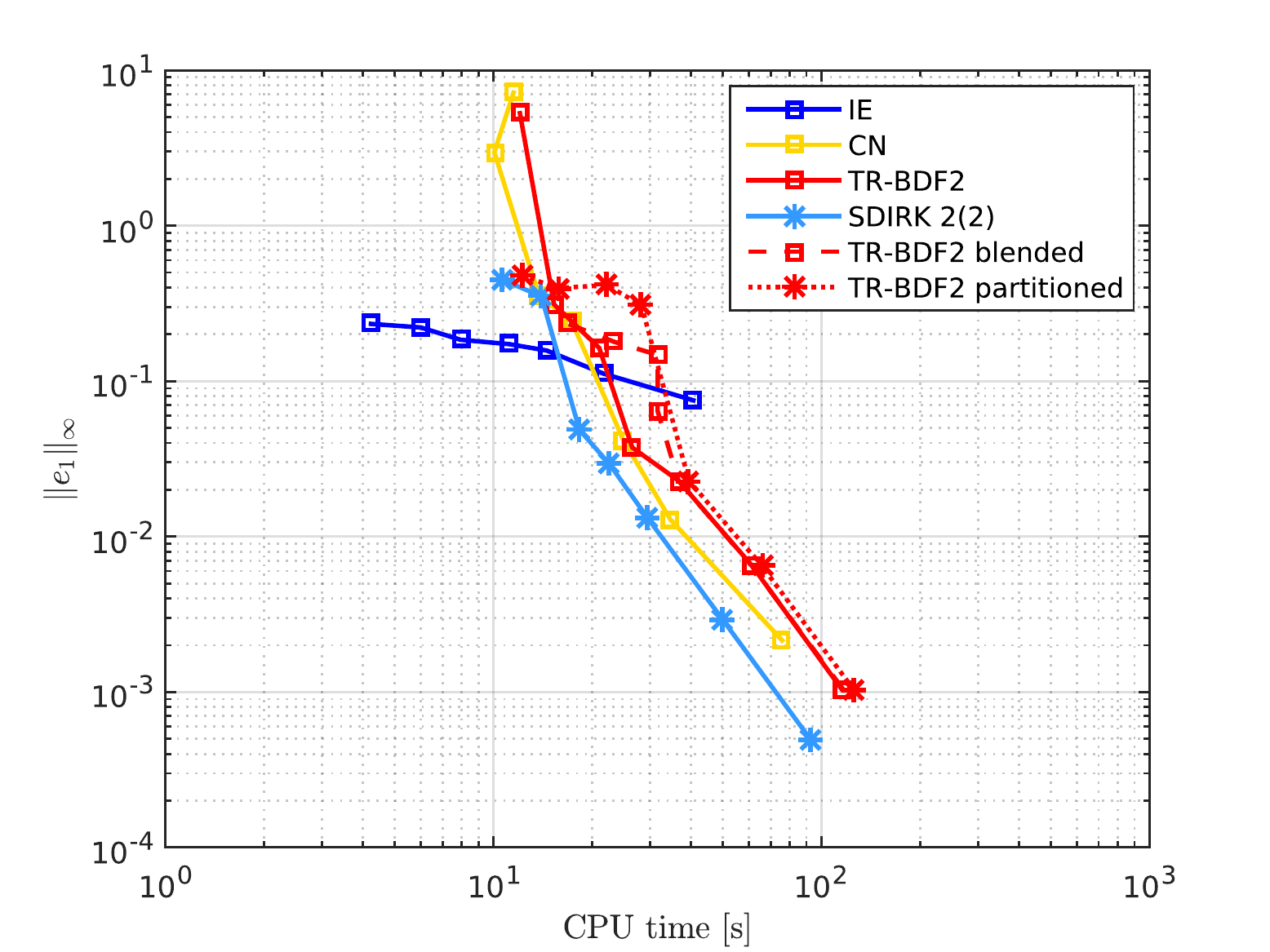}(b)
 \caption{Buckley Leverett equation with Koren limiter and non-smooth initial condition. \emph{Left}:
   Error-stepsize. \emph{Right}: Error-workload.}
 \label{fig:buckley_leverett_inf}
\end{figure}

The accuracy curves for the Buckley-Leverett equation in Figure \ref{fig:buckley_leverett_inf} confirm our
previous findings. Conditionally monotone methods achieve worse results than implicit Euler at coarse step
sizes, due to the impact of the relevant TVD violations. TR-BDF2 blended offers a seamless compromise between
accuracy at fine time steps and monotonicity at coarse time steps, while TR-BDF2 partitioned obtains
qualitatively very good solutions, but it is penalized in the $l^\infty$ norm by the wrong prediction on the
advection speed of the leading shock.

\section{Conclusions}
\label{conclusions}

We have reviewed a general framework for the preservation of some relevant solution properties during
numerical integrations with RK methods. The generality of the absolute monotonicity results proved to be of
practical relevance for assessing contractivity, monotonicity, positivity and strong stability of RK
methods. In particular we analyzed the monotonicity properties of the TR-BDF2 method, that was introduced in
\cite{bank:1985} and successively reformulated and analyzed in \cite{hosea:1996}. We derived the
characteristic SSP coefficient of the DIRK family associated to TR-BDF2, which expresses a CFL-like condition
for monotonicity properties.  We proposed two modifications, the first one based on a hybridization in time
with the implicit Euler method and the second one being an automatically partitioned RK method that tries to
separate monotone and non-monotone solution components in each time step. Both strategies attempt to enforce
monotonicity properties in constant time step integrations, as commonly found in meteorology, environmental
fluid dynamics or turbulent reactive flow simulations. Both monotone strategies make use of sensor functions
able to detect local or global violations of suitable functional bounds, thus triggering a robust integration
procedure when necessary to maintain monotonicity.
 In this way accuracy, is sacrificed locally in order to preserve monotonicity independently of the
time step and stiffness of the problem.

Both strategies were assessed empirically against other RK methods on a series of benchmark problems, ranging
from zero dimensional chemical rectors to advection diffusion reaction equations and nonlinear conservation
laws. The results show that the time hybridization strategies are able to guarantee a seamless compromise between accuracy
at fine step sizes and monotonicity at coarse step sizes, while the partitioned strategy obtains promising
results penalized only by a reduced shock advection speed at high CFL. Further research may be useful to
identify more appropriate sensors for triggering the partitioning method, as well as to extend the two
strategies to SSP-optimal RK methods such as SDIRK 2(2). 

%
%

\section{Appendix}

\subsection{Relevant properties for chemical kinetics problems}

Ordinary differential equations in the form \eqref{eq:ode} arise in the context of chemical kinetics when the
reaction term is modeled by the \emph{mass action law}. By ignoring the presence of additional terms in
advection-diffusion-reaction PDEs, we briefly review here the reaction term.

Considering a chemical system having $N_s$ species interacting in $N_r$ chemical reactions
\begin{equation*}
 \sum\limits_{q=1}^{N_s} l_{pq} u_{q} \xrightarrow[]{k_p} \sum\limits_{q=1}^{N_S} r_{pq} u_{q}
\end{equation*}
the evolution of the molar concentrations is described by the mass action law:
\begin{subequations}
 \label{eq:mass_action_law}
\begin{align}
   \frac{\mathrm{d}}{\mathrm{d}t}u(t) &= \dot \omega \quad \text{and} \quad t\geq t^0, \\
   u(t^0) &= u^0 .
\end{align}
\end{subequations}
The source term is $\dot \omega = Q \omega(u)$ where $Q$ is the matrix of stoichiometric coefficients
\begin{equation*}
 Q = R - L \in \mathbb{R}^{N_s \times N_r} \quad \text{with:} \quad R=\left[ r_{pq} \right] \quad L=\left[ l_{pq}
   \right]
\end{equation*}
and $\omega \in \mathbb{R}^{N_r}$ represents the vector of reaction rates that are usually expressed following the
exponential form of the Arrhenius law
\begin{equation*}
 \omega_p (u) = k_p \prod\limits_{q=1}^{N_s} (u_q)^{l_{pq}} \quad \text{for} \quad p=1,\ldots,N_r.
\end{equation*}
The relevant physical properties for the solution to \eqref{eq:mass_action_law} are:
\begin{itemize}
 \item \emph{Conservation of atomic mass}: the mass of single atomic species such as C, O and N forming the
   chemical species remains constant, since atomic species are conserved during chemical
   reactions. Algebraically this means that if $e \!\in\! \text{Ker}(Q^\intercal)$ is a linear invariant of
   \eqref{eq:mass_action_law} and from $\text{rank}(Q) \!=\! N_s - n$ we know that there are $n$ linear
   invariants, then it is possible to collect the vectors defining the null space of $Q^\intercal$ in the
   columns of the matrix $A \in \mathbb{R}^{N_s \times n}$. As a consequence the solution of
   \eqref{eq:mass_action_law} must satisfy
   \begin{equation}
    \label{eq:mass_conservation}
     A^\intercal u(t) = A^\intercal u^0 = const \quad \text{for} \quad t\geq t^0.
   \end{equation}
 \item \emph{Positivity}: the concentrations are physical quantities that are bounded in the range of
   significant values $0 \!\leq\! u \!\leq\! 1$. By splitting the production and destruction terms in
   \eqref{eq:mass_action_law}
   \begin{subequations}
     \begin{align}
      \frac{\mathrm{d}}{\mathrm{d}t}u(t) &= P(u)-D(u)u \quad \text{and} \quad t\geq t^0, \\
      u(t^0) &= u^0 .
     \end{align}
   \end{subequations}
   where the production and destruction terms are
   \begin{equation*}
    P(u) = R \, \omega(u) \quad \text{and} \quad D(u)=\text{diag}\left[ \frac{L(i,\cdot) \, \omega(u)}{u_p}
     \right].
   \end{equation*}
   This form ensures that $D_{ii}(u)$ are polynomials due to the functional form of the reaction rates. This
   implies that if all concentrations are nonnegative except $u_p \!=\! 0$ then
   \begin{equation*}
    \frac{\mathrm{d}}{\mathrm{d}t}u_p(t) = P_p(u) \geq 0 \implies u(t) \geq 0 \quad \text{whenever}
    \quad u(0) \!\geq\! 0 \, .
   \end{equation*}

\end{itemize}
While conservation of linear invariants is automatic in general linear methods, positivity is more difficult
to achieve. In practical applications positivity is usually enforced by a \emph{clipping} step, in which all
negative solution components are explicitly set to zero. Clipping alters the conservation of mass, since the
error is introduced in a single direction only, namely adding mass to the system. While for short term
computations under tight tolerances this is generally acceptable, for long time integrations the mass added
may be detrimental to solution accuracy. Moreover, as evident from Section~\ref{review}, positivity is more
than a constraint for physical significance of a numerical solution, since it is strongly related to
nonlinear stability properties of numerical methods.

\subsection{Tables of TV seminorm in numerical experiments}

Here below we report the tables of the TV seminorm measured in numerical experiments to show that TV
violations closely follow the expected behaviour from the absolute monotonicity radius.

\begin{table}[htbc]
 \footnotesize
 \centering
 \begin{tabularx}{1.00\textwidth}{c c c c c c}
  \toprule
  $h$ & $Cou$ & $ref$ & $IE$ & $CN$ & $TR\text{-}BDF2 (clip)$ \\
  \midrule
  0.0025 &  0.25 & 2.00000000 & 2.00000000 & 2.00000000 & 2.00000000 \\ 
  0.0050 &  0.50 & 2.00000000 & 2.00000000 & 2.00000000 & 2.00000000 \\ 
  0.0100 &  1.00 & 2.00000000 & 2.00000000 & 2.00000000 & 2.00000000 \\
  0.0200 &  2.00 & 2.00000000 & 2.00000000 & 2.00000000 & 2.00000000 \\
  0.0241 &  2.41 & 2.00000000 & 2.00000000 & 2.37516991 & 2.00000000 \\
  0.0400 &  4.00 & 2.00000000 & 2.00000000 & 3.33333333 & 2.27858017 \\
  0.0600 &  6.00 & 2.00000000 & 2.00000000 & 4.06243821 & 2.39070772 \\
  0.1000 & 10.00 & 2.00000000 & 2.00000000 & 5.21857423 & 2.47739160 \\
  \midrule
  & & $SDIRK 2(2)$ & $ROS2$ & $TR\text{-}BDF2 (blend)$ & $TR\text{-}BDF2 (part.)$ \\
  \midrule
  0.0025 & & 2.00000000 & 2.00877086 & 2.00000000 & 2.00000000 \\
  0.0050 & & 2.00000000 & 2.02925347 & 2.00000000 & 2.00000000 \\
  0.0100 & & 2.00000000 & 2.07630970 & 2.00000000 & 2.00000000 \\
  0.0200 & & 2.00000000 & 2.14215613 & 2.00000000 & 2.00000000 \\
  0.0241 & & 2.00000000 & 2.14775690 & 2.00000000 & 2.00000000 \\
  0.0400 & & 2.00000000 & 2.12378933 & 2.00000000 & 2.00000000 \\
  0.0600 & & 2.76800000 & 2.07354078 & 2.00000000 & 2.00114309 \\
  0.1000 & & 3.73260435 & 2.01991743 & 2.00000000 & 2.00000000 \\
  \bottomrule
 \end{tabularx}
 \caption{$\|TV\|_{\infty}$ for the advection problem with non smooth initial condition.}
 \label{tab:tv_adv_discont_data}
\end{table}
\clearpage
\begin{table}[htbc]
 \footnotesize
 \centering
 \begin{tabularx}{1.00\textwidth}{c c c c c c}
  \toprule
  $h$ & $Cou$ & $ref$ & $IE$ & $CN$ & $TR\text{-}BDF2 (clip)$ \\
  \midrule
  0.0025 & 0.025 & 19.96000000 & 19.96000000 & 19.96000000 & 19.96000000 \\ 
  0.0050 & 0.050 & 19.96000000 & 19.96000000 & 19.96000000 & 19.96000000 \\ 
  0.0100 & 0.100 & 19.96000000 & 19.96000000 & 19.96000000 & 19.96000000 \\
  0.0250 & 0.200 & 19.96000000 & 19.96000000 & 19.96000000 & 19.96000000 \\
  0.0500 & 0.241 & 19.96000000 & 19.96000000 & 19.96000000 & 19.96000000 \\
  0.1000 & 1.000 & 19.96000000 & 19.96000000 & 21.26167041 & 19.96000000 \\
  \midrule
  & & $SDIRK 2(2)$ & $ROS2$ & $TR\text{-}BDF2 (blend)$ & $TR\text{-}BDF2 (part.)$ \\
  \midrule
  0.0025 & & 19.96000000 & 19.96178624 & 19.96000000 & 19.96000000 \\
  0.0050 & & 19.96000000 & 19.96634023 & 19.96000000 & 19.96000000 \\
  0.0100 & & 19.96000000 & 19.97806568 & 19.96000000 & 19.96000000 \\
  0.0250 & & 19.96000000 & 20.01749991 & 19.96000000 & 19.96000000 \\
  0.0500 & & 19.96000000 & 20.07270798 & 19.96000000 & 19.96000000 \\
  0.1000 & & 19.96000000 & 20.11087432 & 19.96000000 & 19.96000000 \\
  \bottomrule
 \end{tabularx}
 \caption{$\|TV\|_{\infty}$ for the advection diffusion reaction problem with non smooth initial condition.}
 \label{tab:tv_adv_diff_react_discont_data}
\end{table}

\begin{table}[!htbp]
 \footnotesize
 \centering
 \begin{tabularx}{1.00\textwidth}{c c c c c c}
  \toprule
  $h$ & $Cou_{max}$ & $ref$ & $IE$ & $CN$ & $TR\text{-}BDF2 (clip)$ \\
  \midrule
  0.0025 & 0.188 & 1.00000000 & 1.00000000 & 1.00000000 & 1.00000000 \\
  0.0050 & 0.375 & 1.00000000 & 1.00000000 & 1.00000000 & 1.00000000 \\
  0.0100 & 0.750 & 1.00000000 & 1.00000000 & 1.00000000 & 1.00000000 \\
  0.0200 & 1.500 & 1.00000000 & 1.00000000 & 1.06702208 & 1.00000000 \\
  0.0400 & 3.000 & 1.00000000 & 1.00000000 & 1.36710392 & 1.10154954 \\
  0.0600 & 4.500 & 1.00000000 & 1.00000000 & 1.33460354 & 1.18618300 \\
  0.1000 & 7.500 & 1.00000000 & 1.00000000 & 1.20438152 & 1.17675243 \\
  \midrule
  & & $EE$ & $SDIRK 2(2)$ & $TR\text{-}BDF2 (blend)$ & $TR\text{-}BDF2 (part.)$\\
  \midrule
  0.0025 & & 1.00000000 & 1.00000000 & 1.00000000 & 1.00000000 \\
  0.0050 & & 1.00000000 & 1.00000000 & 1.00000000 & 1.00000000 \\
  0.0100 & & 3.47978819 & 1.00000000 & 1.00000000 & 1.00000000 \\
  0.0200 & & $\infty$   & 1.00000000 & 1.00000000 & 1.00000000 \\
  0.0400 & & $\infty$   & 1.02205723 & 1.00000000 & 1.00000000 \\
  0.0600 & & $\infty$   & 1.39639398 & 1.00000000 & 1.00000000 \\
  0.1000 & & $\infty$   & 1.85272884 & 1.00000000 & 1.00000000 \\
  \bottomrule
 \end{tabularx}
 \caption{$\|TV\|_{\infty}$ for the Burgers equation with van Leer limiter and smooth initial condition.}
 \label{tab:tv_burgers_data}
\end{table}

\begin{table}[!htbp]
 \footnotesize
 \centering
 \begin{tabularx}{1.00\textwidth}{c c c c c c}
  \toprule
  $h$ & $Cou_{max}$ & $ref$ & $IE$ & $CN$ & $TR\text{-}BDF2 (clip)$ \\
  \midrule
  0.0010 & 0.2206 & 1.00000000 & 1.00000000 &  1.00000000 &  1.00000000 \\
  0.0025 & 0.5514 & 1.00000000 & 1.00000000 &  1.00000000 &  1.00000000 \\
  0.0050 & 1.1007 & 1.00000000 & 1.00000001 &  1.00000000 &  1.00000000 \\
  0.0075 & 1.6543 & 1.00000000 & 1.00000000 &  1.27084161 &  1.04519501 \\
  0.0100 & 2.2007 & 1.00000000 & 1.00000000 &  1.40310073 &  1.18671375 \\
  0.0150 & 3.3010 & 1.00000000 & 1.00000000 &  8.09702136 &  1.32271889 \\
  0.0250 & 5.5037 & 1.00000000 & 1.00000000 & 16.39900045 & 12.48391471 \\
  \midrule
  & & $EE$ & $SDIRK 2(2)$ & $TR\text{-}BDF2 (blend)$ & $TR\text{-}BDF2 (part.)$ \\
  \midrule
  0.0010 & &  1.00000000 & 1.00000000 & 1.00000000 &  1.00000000 \\
  0.0025 & &  1.00000000 & 1.00000000 & 1.00000000 &  1.00000000 \\
  0.0050 & &  1.87710600 & 1.00000000 & 1.00000000 &  1.00000000 \\
  0.0075 & &  6.74722895 & 1.00000000 & 1.00000000 &  1.00000000 \\
  0.0100 & & 21.35210300 & 1.00000000 & 1.00000001 &  1.00495631 \\
  0.0150 & & 21.10782320 & 1.28222058 & 1.00000000 &  1.03092403 \\
  0.0250 & & 19.65284356 & 1.57372343 & 1.00000000 &  1.00663983 \\
  \bottomrule
 \end{tabularx}
 \caption{$\|TV\|_{\infty}$ for the Buckley-Leverett equation with Koren limiter and non smooth initial
   condition.}
 \label{tab:tv_buckley_leverett_data}
\end{table}

\section*{Acknowledgements}

A.D.R. would like to thank Tenova S.p.A. for sponsoring his Executive PhD at Politecnico di Milano and all the
faculty members at the Department of Mathematics for their continuous support.  L.B. acknowledges financial support
from the INDAM -  GNCS project 'Metodi numerici semi-impliciti e semi-Lagrangiani per sistemi iperbolici di leggi di bilancio'.
Useful discussions with L. Formaggia and A. Scotti on the topics studied in this paper are kindly acknowledged.

\bibliographystyle{plain}
\bibliography{SSPbibliography}

\end{document}